\documentclass[11pt]{article}
\usepackage{amssymb}
\usepackage{latexsym}
\usepackage{amsfonts}
\usepackage{amsmath}
\newtheorem{thrm}{Theorem}[section]
\newtheorem{lemma}[thrm]{Lemma}
\newtheorem{prop}[thrm]{Proposition}

\newtheorem{cor}[thrm]{Corollary}
\newtheorem{remark}[thrm]{Remark}

\numberwithin{equation}{section}

\usepackage[dvips]{color}

\def\P{\mathbb{P} }

\def\R{\mathbb{R} }
\def\V{\mathbb{V} }
\def\N{\mathbb{N} }
\def\D{\mathbb{D} }
\def\C{\mathcal{C} }

\begin{document}
\allowdisplaybreaks
\begin{doublespace}
\title{Functional central limit theorems for supercritical superprocesses}
\author{ \bf  Yan-Xia Ren\footnote{The research of this author is supported by NSFC (Grant No.  11271030 and 11128101) and Specialized Research Fund for the
Doctoral Program of Higher Education.\hspace{1mm} } \hspace{1mm}\hspace{1mm}
Renming Song\thanks{Research supported in part by a grant from the Simons
Foundation (208236).} \hspace{1mm}\hspace{1mm} and \hspace{1mm}\hspace{1mm}
Rui Zhang
\hspace{1mm} }
\date{}
\maketitle
\begin{abstract}
In this paper, we establish some functional central limit theorems for a large class of general supercritical
superprocesses with spatially dependent branching mechanisms satisfying a second moment condition.
In the particular case when the state $E$ is a finite set and the underline motion is an
irreducible Markov chain on $E$, our results are superprocess analogs of
 the functional central limit theorems of \cite{Janson} for supercritical multitype branching processes.
The results of this paper are
refinements of the central limit theorems in \cite{RSZ3}.
\end{abstract}

\medskip
\noindent {\bf AMS Subject Classifications (2000)}: Primary 60J68;
Secondary 60F05, 60G57, 60J45

\medskip

\noindent{\bf Keywords and Phrases}: Functional central limit theorem,
supercritical superprocess, excursion measures of superprocesses.

\bigskip

\baselineskip=6.0mm

\section{Introduction}

Kesten and Stigum \cite{KS, KS66} initiated the study of central limit theorems for supercritical branching processes.
In these two papers, they established central limit theorems for
supercritical multi-type Galton-Watson processes by using the
Jordan canonical form of the mean matrix.
Then in \cite{Ath69a, Ath69, Ath71}, Athreya proved central limit
theorems for supercritical multi-type continuous time branching processes,
also using the Jordan canonical form of the mean matrix.
Asmussen and Keiding \cite{AK} used martingale central limit theorems to prove central limit theorems for
supercritical multi-type branching processes. In \cite{AH83}, Asmussen and Hering established
spatial central limit theorems for general supercritical branching Markov processes under a certain
condition. In \cite{Janson}, Janson extended the results of \cite{Ath69a, Ath69, Ath71, KS, KS66}
and established functional central limit theorems for  multitype branching processes.
In \cite[Remark 4.1]{Janson}, Janson mentioned the possibility of extending his functional
central limit theorems to the case of infinitely many types (with suitable assumptions). However,
he ended this remark with the following sentence:
``It is far from clear how such an extension should be formulated, and
we have not pursued this''.

The recent study of spatial central limit theorems for branching Markov processes started with
\cite{RP}.
In this paper, Adamczak and Mi{\l}o\'{s} proved some central limit theorems for supercritical branching
Ornstein-Uhlenbeck processes with binary branching mechanism.
In \cite{Mi}, Mi{\l}o\'{s} proved
some central limit theorems for supercritical super Ornstein-Uhlenbeck processes with branching
mechanisms satisfying a fourth moment condition.
In \cite{RSZ}, we established central limit
theorems for supercritical super Ornstein-Uhlenbeck processes with
branching mechanisms satisfying only a second moment condition. More importantly,
compared with the results of
\cite{RP, Mi}, the central limit
theorems in \cite{RSZ} are more satisfactory since our limit normal random variables are non-degenerate.
In \cite{RSZ2}, we sharpened and generalized the spatial central limit
theorems mentioned above, and obtained central limit theorems for a large class
of general supercritical branching symmetric Markov processes with  spatially dependent branching mechanisms satisfying
only a second moment condition. In \cite{RSZ3}, we obtained central limit theorems for a large class of
general supercritical superprocesses with symmetric spatial motions and with spatially dependent branching mechanisms satisfying
only a second moment condition. Furthermore, we also obtained the covariance structure of the limit
Gaussian field in \cite{RSZ3}. In \cite{RSZ4}, we extended the results of \cite{RSZ2} to
supercritical branching nonsymmetric Markov processes with
spatially dependent branching mechanisms satisfying
only a second moment condition.

The main purpose of this  paper is
to establish functional central limit theorems, for supercritical superprocesses
with spatially dependent branching mechanisms satisfying
only a second moment condition, similar to those
of \cite{Janson}, for supercritical multitype branching processes.
For simplicity, we will assume the
spatial process is symmetric. One could combine the techniques of this
paper with that of \cite{RSZ4} to extend the results of this paper to the case when the spatial motion is
not symmetric. We leave this to the interested reader.

The organization of this paper is as follows.
In the remainder of this section, we spell out our assumptions and present
our main result.
Section 2 contains some preliminary results, while the
proof of the main result is given in Section 3.

\subsection{Spatial process}\label{subs:sp}

Our assumptions on the underlying spatial process are the same as in \cite{RSZ2}.
In this subsection, we recall the assumptions on the spatial process.

$E$ is a locally compact separable metric space and $m$ is a $\sigma$-finite Borel measure
on $E$ with full support.
$\partial$ is a point not contained in $E$ and will be interpreted as the cemetery point.
Every function $f$ on $E$ is automatically extended to $E_{\partial}:=E\cup\{\partial\}$
by setting $f(\partial)=0$.
We will assume that $\xi=\{\xi_t,\Pi_x\}$ is an $m$-symmetric
Hunt process on $E$.
The semigroup of $\xi$ will be denoted by $\{P_t:t\geq 0\}$.
We will always assume
that there exists a family of
continuous strictly positive symmetric functions $\{p_t(x,y):t>0\}$ on $E\times E$ such that
$$
  P_tf(x)=\int_E p_t(x,y)f(y)\,m(dy).
$$
It is well-known that for $p\geq 1$, $\{P_t:t\ge 0\}$ is a
strongly continuous contraction semigroup on $L^p(E, m)$.

Define $\widetilde{a}_t(x):=p_t(x, x)$.
We will always assume that $\widetilde{a}_t(x)$ satisfies the following two conditions:
\begin{description}
  \item[(a)] For any $t>0$, we have
$$
        \int_E \widetilde{a}_t(x)\,m(dx)<\infty.
$$
   \item[(b)] There exists $t_0>0$ such that $\widetilde{a}_{t_0}(x)\in L^2(E,\,m)$.
\end{description}
It is easy to check (see \cite{RSZ2}) that condition $(b)$ above is equivalent to
 \begin{description}
   \item[(b$'$)] There exists $t_0>0$ such that for all $t\ge t_0$, $\widetilde{a}_{t}(x)\in L^2(E,m)$.
 \end{description}

These two conditions are satisfied by a lot of Markov processes. In \cite{RSZ2}, we gave several classes
of examples of Markov processes satisfying these two conditions.

\subsection{Superprocesses}

Our basic assumptions on the superprocess are the same as in \cite{RSZ3}.
In this subsection, we recall these assumptions.
Let $\mathcal{B}_b(E)$ ($\mathcal{B}_b^+(E)$) be the set of
(nonnegative) bounded  Borel functions on $E$.

The superprocess $X=\{X_t:t\ge 0\}$
is determined by three parameters:
a spatial motion $\xi=\{\xi_t, \Pi_x\}$ on $E$ satisfying the
assumptions of the previous subsection,
a branching rate function $\beta(x)$ on $E$ which is a nonnegative bounded Borel
function and a branching mechanism $\psi$ of the form
\begin{equation}
\psi(x,\lambda)=-a(x)\lambda+b(x)\lambda^2+\int_{(0,+\infty)}(e^{-\lambda y}-1+\lambda y)n(x,dy),
\quad x\in E, \quad\lambda> 0,
\end{equation}
where $a\in \mathcal{B}_b(E)$, $b\in \mathcal{B}_b^+(E)$ and $n$ is a kernel from $E$ to $(0,\infty)$ satisfying
\begin{equation}\label{n:condition}
  \sup_{x\in E}\int_0^\infty y^2 n(x,dy)<\infty.
\end{equation}

Let ${\cal M}_F(E)$
be the space of finite measures on $E$, equipped with topology of weak convergence.
The superprocess $X$ is a  Markov process taking values in ${\cal M}_F(E)$.
The existence of such superprocesses is well-known, see, for instance,
\cite{E.B.} or \cite{Li11}.
As usual, $\langle f,\mu\rangle:=\int f(x)\mu(dx)$ and $\|\mu\|:=\langle 1,\mu\rangle$.
According to \cite[Theorem 5.12]{Li11}, there is a Borel right process $X=\{\Omega, {\cal G}, {\cal G}_t, X_t, \P_\mu\}$ taking values in  $\mathcal{M}_F(E)$ such that
for every
$f\in \mathcal{B}^+_b(E)$ and $\mu \in \mathcal{M}_F(E)$,
\begin{equation}
  -\log \P_\mu\left(e^{-\langle f,X_t\rangle}\right)=\langle u_f(\cdot,t),\mu\rangle,
\end{equation}
where $u_f(x,t)$ is the unique positive solution to the equation
\begin{equation}\label{1.3}
  u_f(x,t)+\Pi_x\int_0^t\psi(\xi_s, u_f(\xi_s,t-s))\beta(\xi_s)ds=\Pi_x f(\xi_t),
\end{equation}
where $\psi(\partial,\lambda)=0, \lambda>0$.
By the definition of Borel right processes (see \cite[Definition A.18]{Li11}),
$({\cal G}, {\cal G}_t)_{t\ge 0}$ are augmented, $(\mathcal{G}_t:t\ge 0)$ is right continuous
and $X$ satisfies the Markov property with respect to $(\mathcal{G}_t:t\ge 0)$.
Moreover, such a superprocess $X$ has a Hunt realization in ${\cal M}_F(E)$,
see \cite[Theorem 5.12]{Li11}.
In this paper, the superprocess we deal with is always this Hunt realization.

Define
\begin{equation}\label{e:alpha}
\alpha(x):=\beta(x)a(x)\quad \mbox{and }
A(x):=\beta(x)\left( 2b(x)+\int_0^\infty y^2 n(x,dy)\right).
\end{equation}
Then, by our assumptions, $\alpha(x)\in \mathcal{B}_b(E)$ and $A(x)\in\mathcal{B}_b(E)$.
Thus there exists $K>0$ such that
\begin{equation}\label{1.5}
  \sup_{x\in E}\left(|\alpha(x)|+A(x)\right)\le K.
\end{equation}
For any $f\in\mathcal{B}_b(E)$ and $(t, x)\in (0, \infty)\times E$, define
\begin{equation}\label{1.26}
   T_tf(x):=\Pi_x \left[e^{\int_0^t\alpha(\xi_s)\,ds}f(\xi_t)\right].
\end{equation}
It is well-known that $T_tf(x)=\P_{\delta_x}\langle f,X_t\rangle$ for every $x\in E$.

It is shown in \cite{RSZ2}
that there exists a family of continuous strictly positive symmetric functions
$\{q_t(x,y),t>0\}$ on $E\times E$ such that
$q_t(x,y)\le e^{Kt}p_t(x,y)$ and
for any $f\in \mathcal{B}_b(E)$,
$$
  T_tf(x)=\int_E q_t(x,y)f(y)\,m(dy).
$$
It follows immediately that, for any $p\ge 1$, $\{T_t: t\ge 0\}$ is a
strongly continuous semigroup on $L^p(E, m)$ and
\begin{equation}\label{Lp}
  \|T_tf\|_p^p\le e^{pKt}\|f\|_p^p.
\end{equation}

Define $a_t(x):=q_t(x, x)$. It follows from the assumptions (a) and (b) in the
previous subsection that $a_t$ enjoys the following properties:
\begin{description}
  \item[(i)] For any $t>0$, we have
$$
        \int_E a_t(x)\,m(dx)<\infty.
$$
  \item[(ii)] There exists $t_0>0$ such that for all $t\ge t_0$, $a_{t}(x)\in L^2(E,m)$.
\end{description}

By H\"{o}lder's inequality, we get
$$q_t(x,y)=\int_E q_{t/2}(x,z)q_{t/2}(z,y)\,m(dz)\le a_t(x)^{1/2}a_t(y)^{1/2}.$$
Since $q_t(x,y)$ and $a_t(x)$ are continuous in $x\in E$,
by the dominated convergence theorem, we get that, if $f\in L^2(E,m)$,
$T_tf(\cdot)$ is continuous for any $t>0$.

It follows from (i) above that, for any $t>0$, $T_t$ is  a compact
operator. The infinitesimal generator
$\mathcal{L}$ of $\{T_t:t\geq 0\}$ in $L^2(E, m)$ has purely discrete spectrum with eigenvalues
$-\lambda_1>-\lambda_2>-\lambda_3>\cdots$.
It is known that either the number of these eigenvalues is finite, or $\lim_{k\to\infty}\lambda_k=\infty$.
The first eigenvalue $-\lambda_1$ is simple and the eigenfunction
$\phi_1$ associated with $-\lambda_1$ can be chosen to be strictly positive everywhere and continuous.
We will assume that $\|\phi_1\|_2=1$. $\phi_1$ is sometimes denoted as $\phi^{(1)}_1$.
For $k>1$, let $\{\phi^{(k)}_j,j=1,2,\cdots n_k\}$ be
an orthonormal basis of the eigenspace
associated with $-\lambda_k$.
It is well-known that $\{\phi^{(k)}_j,j=1,2,\cdots n_k; k=1,2,\dots\}$ forms a complete orthonormal basis of $L^2(E,m)$
and all the eigenfunctions are continuous.
For any $k\ge 1$, $j=1, \dots, n_k$ and $t>0$, we have $T_t\phi^{(k)}_j(x)=e^{-\lambda_k t}\phi^{(k)}_j(x)$ and
\begin{equation}
\label{1.37}
e^{-\lambda_kt/2}|\phi^{(k)}_j|(x)\le a_t(x)^{1/2}, \qquad x\in E.
\end{equation}
It follows from the relation above that all the eigenfunctions $\phi^{(k)}_j$ belong to $L^4(E, m)$.
The basic facts recalled in this paragraph are well-known, for instance, one can
refer to  \cite[Section 2]{DS}.

In this paper, we always assume that the superprocess $X$ is supercritical,
that is, $\lambda_1<0$.

In this paper, we also assume that, for any $t>0$ and $x\in E$,
\begin{equation}\label{extinction}
  \P_{\delta_x}\{\|X_t\|=0\}\in(0,1).
\end{equation}
Here is a sufficient condition for \eqref{extinction}.
Suppose that $\Phi(z)=\inf_{x\in E}\psi(x,z)\beta(x)$ can be written in the form:
$$
\Phi(z)=\widetilde{a}z+\widetilde{b}z^2+\int^\infty_0(e^{-zy}-1+zy)\widetilde{n}(dy)
$$
with $\widetilde{a}\in \R$, $\widetilde{b}\ge 0$ and $ \widetilde{n}$ being a measure on $(0,\infty)$ satisfying
$\int^\infty_0(y\wedge y^2)\widetilde{n}(dy)<\infty$.
If $\widetilde{b}+\widetilde{n}(0,\infty)>0$ and $\Phi(z)$ satisfies
\begin{equation}\label{Phi}
  \int^\infty\frac{1}{\Phi(z)}\,dz<\infty,
\end{equation}
then \eqref{extinction} holds. For the last claim, see, for instance, \cite[Lemma 11.5.1]{Dawson}.

\subsection{Main Result}

In the remainder of this paper,
whenever we deal with
an initial configuration
$\mu\in {\cal M}_F(E)$,
we are implicitly assuming that it has compact support.

We will use $(\cdot, \cdot)_m$ to denote inner product in
$L^2(E, m)$.
Any $f\in L^2(E,m)$ admits the following expansion:
$$
  f(x)=\sum_{k=1}^\infty\sum^{n_k}_{j=1} a_j^k\phi_j^{(k)}(x),
$$
where $a_j^k=(f,\phi_j^{(k)})_m$ and the series converges in $L^2(E,m)$.
$a^1_1$ will sometimes be written as $a_1$.
For $f\in L^2(E,m)$, define
$$
  \gamma(f):=\inf\{k\geq 1: \mbox{ there exists } j \mbox{ with }
  1\leq j\leq n_k\mbox{ such that }
  a_j^k\neq 0\},
$$
where we use the usual convention $\inf\varnothing=\infty$.
We note that if $f\in L^2(E,m)$ is nonnegative and $m(x: f(x)>0)>0$,
then $(f,\phi_1)_m>0$,
which implies $\gamma(f)=1$.

Define
\begin{equation*}
  H_t^{k,j}:=e^{\lambda_k t}\langle\phi_j^{(k)}, X_t\rangle,\quad t\geq 0.
\end{equation*}
In \cite[Lemma 1.1]{RSZ3}, it has been proved that, for any nonzero $\mu\in {\cal M}_F(E)$,  $H_t^{k,j}$ is a martingale under $\P_{\mu}$.
Moreover, if $\lambda_1>2\lambda_k$, then $\sup_{t>3t_0}\P_{\mu}(H_t^{k,j})^2<\infty$.
Thus the limit
\begin{equation*}
  H_\infty^{k,j}:=\lim_{t\to \infty}H_t^{k.j}
\end{equation*}
exists $\P_{\mu}$-a.s. and in $L^2(\P_{\mu})$.

In particular, we write $W_t:=H^{1,1}_t=e^{\lambda_1 t}\langle \phi_1, X_t\rangle$ and $W_\infty:=H^{1,1}_\infty$.
$\{W_t: t\ge 0\}$ is a nonnegative martingale and
$$
  W_t \to W_\infty,\quad \P_{\mu}\mbox{-a.s. and in }L^2(\P_{\mu}).
$$
Thus $W_\infty$ is non-degenerate. Moreover, we have $\P_{\mu}(W_\infty)=\langle\phi_1, \mu\rangle$.
Put $\mathcal{E}=\{W_\infty=0\}$, then $\P_{\mu}(\mathcal{E})<1$.
It is clear that  $\mathcal{E}^c\subset\{X_t(E)>0,\forall t\ge 0\}$.

The following three subspaces of $L^2(E, m)$ will be needed in the statement of the main result:
$$
\C_l:=\left\{g(x)=\sum_{k:\lambda_1>2\lambda_k}\sum^{n_k}_{j=1} b_j^k\phi_j^{(k)}(x): b_j^k\in \R\right\},
$$
$$
\C_c:=\left\{g(x)=\sum^{n_k}_{j=1} b_j^k\phi_j^{(k)}(x): 2\lambda_k=\lambda_1, b_j^k\in \R\right\}
$$
and
$$
\C_s:=\left\{g(x)\in L^2(E,m)\cap L^4(E,m):\lambda_1<2\lambda_{\gamma(g)}\right\}.
$$
The space $\C_l$ consists of the functions in $L^2(E, m)$ that only have nontrivial
projections onto the eigen-spaces corresponding to those ``large'' eigenvalues $-\lambda_k$ satisfying
$\lambda_1>2\lambda_k$. The space $\C_l$ is of finite dimension. The space
$\C_c$ is the (finite dimensional) eigen-space corresponding to
the ``critical'' eigenvalue $-\lambda_k$ with  $\lambda_1=2\lambda_k$. Note that there may not
be a critical eigenvalue and $\C_c$ is empty in this case. The space
$\C_s$ consists of the functions in $L^2(E, m)\cap L^4(E, m)$ that only have nontrivial
projections onto the eigen-spaces corresponding to those ``small'' eigenvalues $-\lambda_k$ satisfying
$\lambda_1<2\lambda_k$. The space $\C_s$ is of infinite dimensional in general.

Fix a $q>\max\{K, -2\lambda_1\}$.
For any $p\ge 1$ and  $f\in L^p(E,m)$, define
$$U_q|f|(x):=\int_0^\infty e^{-qs}T_s(|f|)(x)\,ds,\quad x\in E.$$
Then,
\begin{equation}\label{L_p}
  \left(\int_E (U_q|f|(x))^p\,m(dx)\right)^{1/p}\le \int_0^\infty e^{-qs}\|T_s(|f|)\|_p\,ds\le \int_0^\infty e^{-qs}e^{Ks}\,ds\|f\|_p<\infty,
\end{equation}
which implies that $U_q|f|\in L^p(E,m)$.
Let $f^+$ and $f^-$ be the positive part and negative part of $f$ respectively.
For any $x\in E$ with $U_q|f|(x)<\infty$, we define
$$U_qf(x):=\int_0^\infty e^{-qs}T_sf(x)\,ds=U_q(f^+)(x)-U_q(f^-)(x),$$
otherwise we define $U_qf(x)$ be an arbitrary real number.
It follows from \eqref{L_p} that $U_q$ is a bounded linear operator on $L^p(E,m)$.
Notice that
$$U_q(\phi_j^{(k)})(x)=(q+\lambda_k)^{-1}\phi_j^{(k)}(x).$$
One can easily check that, for $f\in L^2(E,m)$, $\gamma(U_qf)=\gamma(f)$.
In fact, by Fubini's theorem, we have
\begin{equation}\label{7.41}
 (U_qf,\phi_j^{(k)})_m=\int_0^\infty e^{-qu}(T_uf,\phi_j^{(k)})_m\,du
  =(q+\lambda_k)^{-1}(f,\phi_j^{(k)})_m.
\end{equation}

For any $f\in L^2(E,m)$, the random variable $\langle U_q|f|,X_t\rangle\in [0,\infty]$ is well defined.
Since $\mu$ has compact support and $T_t(U_q|f|)$ is continuous,
$\P_{\mu}(\langle U_q|f|,X_t\rangle)=\langle T_t(U_q|f|),\mu\rangle<\infty$, and thus
$\P_{\mu}(\langle U_q|f|,X_t\rangle<\infty)=1$.
Therefore, for $t\ge 0$,
$\P_{\mu}\left(\langle U_qf,X_t\rangle \mbox{ is finite}\right)=1.$
In Subsection \ref{subs:cad}, we will give a stronger result:
for any $\mu\in {\cal M}_F(E)$ and $f\in L^2(E,m)$, it holds that
$$\P_{\mu}\left(\langle U_q|f|,X_t\rangle<\infty, \forall t\ge0\right
)=\P_{\mu}\left(\langle U_qf,X_t\rangle \mbox{ is finite}, \forall t\ge0\right)=1.$$

We denote by $\D(\R^d)$ the space of all cadlag functions from $[0, \infty)$ into $\R^d$, equipped with the Skorokhod topology.
There is a metric $\delta$ on $\D(\R^d)$ which is compatible with the Skorokhod topology.
See, for instance, \cite[Chapter VI, 1.26]{J.J.}, for the definition of $\delta$.
In the present paper, we will consider weak convergence of processes in the Skorokhod space $\D(\R^d)$,
which is stronger than convergence in finite dimensional distributions.

\smallskip

For $f\in \C_s$, define
\begin{equation}\label{3.10}
  \sigma_{f,\tau}:=e^{\lambda_1\tau/2}\int_{0}^\infty e^{\lambda_1s}
   (A(T_sf)(T_{s+\tau}f),\phi_1)_m\,ds.
\end{equation}
We write $\sigma_{f,0}$ as $\sigma_f^2$.
For $h\in \C_c$, define
\begin{equation}\label{e:rho}
\rho_h^2:=(Ah^2,\phi_1)_m.
\end{equation}
For $g(x)=\sum_{k: 2\lambda_k<\lambda_1}\sum_{j=1}^{n_k}b_j^k\phi_j^{(k)}(x)\in \C_l$, we put
$$ I_ug(x):=\sum_{k: 2\lambda_k<\lambda_1}\sum_{j=1}^{n_k}e^{\lambda_ku}b_j^k\phi_j^{(k)}(x),\quad x\in E,  u\ge 0,$$
and
$$
F_t(g):=\sum_{k:2\lambda_k<\lambda_1}\sum_{j=1}^{n_k}e^{-\lambda_kt}b_j^kH_\infty^{k,j},\quad t\ge 0.
$$
Define
\begin{equation}\label{e:beta}
\beta_{g,\tau}:=e^{-\lambda_1\tau/2}\int_0^\infty e^{-\lambda_1s}
(A(I_sg)(I_{s+\tau}g),\phi_1)_m\,ds.
\end{equation}
We write $\beta_g^2:=\beta_{g,0}$.
For $f\in\C_s$ and $g\in\C_l$, we define
\begin{equation}\label{eta}
  \eta_{\tau_1,\tau_2}(f,g):=-e^{\lambda_1(\tau_1+\tau_2)/2}\int_{\tau_1}^{\tau_2}e^{-\lambda_1u}\big(A(T_{\tau_2-u}f)(I_{u-\tau_1}g),\phi_1\big)_m\,du,\quad 0\le \tau_1\le \tau_2.
\end{equation}

\begin{thrm}\label{Thm1}
Assume that $f\in\C_s$,
$h\in\C_c$, $g\in\C_l$  and $\mu\in \mathcal{M}_F(E)$.
For any $t>0$, define
$$
Y^{1,f}_t(\tau):=e^{\lambda_1(t+\tau)/2}\langle f,X_{t+\tau}\rangle, \quad\tau\ge 0,
$$
$$
Y^{2,h}_t(\tau):=t^{-1/2}e^{\lambda_1(t+\tau)/2}\langle h,X_{t+\tau}\rangle,\quad \tau\ge 0,
$$
and
$$
Y^{3,g}_t(\tau):=e^{\lambda_1(t+\tau)/2}\left(\langle g,X_{t+\tau}-F_{t+\tau}(g)\rangle\right),\quad \tau\ge 0.
$$
Then, for each fixed $t\in [0, \infty)$, $\left(W_t, Y^{1,U_qf}_t(\cdot),Y^{2,h}_t(\cdot),Y^{3,g}_t(\cdot)\right)$
is a $\mathbb{D}(\R^4)$-valued random variable under $\P_\mu$,
where $W_t$ is regarded as a constant process.
Furthermore, under $\P_{\mu}$,
\begin{equation}\label{result}
   \left(W_t, Y^{1,U_qf}_t(\cdot),Y^{2,h}_t(\cdot),Y^{3,g}_t(\cdot)\right)\stackrel{d}{\rightarrow}\left(W_\infty, \sqrt{W_\infty}G^{1,U_qf}(\cdot),\sqrt{W_\infty}G^{2,h},\sqrt{W_\infty}G^{3,g}(\cdot)\right),
      \quad \mbox{ as } t\to\infty,
\end{equation}
in $\D(\R^4)$.
 Here $G^{2,h}\sim \mathcal{N}(0,\rho^2_h)$ is a constant
 process, and $\{(G^{1,U_qf}(\tau),G^{3,g}(\tau)):\tau\ge0\}$ is
 a continuous $\R^2$-valued Gaussian process with mean $0$ and covariance
 functions given by
\begin{equation}\label{3.3}
 E(G^{1,U_qf}(\tau_1)G^{1,U_qf}(\tau_2))=\sigma_{U_qf,\tau_2-\tau_1},\quad \mbox{for }0\le \tau_1\le \tau_2,
\end{equation}
\begin{equation}\label{3.4}
  E(G^{3,g}(\tau_1)G^{3,g}(\tau_2))=\beta_{g,\tau_2-\tau_1},\quad \mbox{for }0\le \tau_1\le \tau_2,
\end{equation}
and
\begin{equation}\label{3.6}
  E(G^{3,g}(\tau_1)G^{1,U_qf}(\tau_2))=\left\{\begin{array}{ll}\displaystyle\eta_{\tau_1,\tau_2}(U_qf,g),&\quad\mbox{if }0\le \tau_1< \tau_2,\\
  \displaystyle 0, &\quad \mbox{if }\tau_1\ge\tau_2\ge 0.\end{array}\right.
\end{equation}
Moreover, $W_\infty$, $G^{2,h}$ and
$(G^{1,U_qf},G^{3,g})$ are independent.
\end{thrm}

\bigskip

For $f\in L^2(E,m)$, we define
\begin{eqnarray*}
     f_{(s)}(x)& :=& \sum_{k:\lambda_1>2\lambda_k}
  \sum_{j=1}^{n_k}a_j^k\phi_j^{(k)}(x),\\
  f_{(l)}(x)&:=& \sum_{k:\lambda_1<2\lambda_k}
   \sum_{j=1}^{n_k}a_j^k\phi_j^{(k)}(x),\\
  f_{(c)}(x)&:=&f(x)- f_{(s)}(x)-f_{(l)}(x).
\end{eqnarray*}
Then $f_{(l)}\in\C_{s}$, $f_{(c)}\in\C_{c}$ and $f_{(s)}\in\C_{l}$.

\begin{remark}\label{large}
Assume that $g=U_qf$ for some $f\in L^2(E,m)\cap L^4(E,m)$ satisfying
$\lambda_1\ge 2\lambda_{\gamma(f)}$.
Then $g_{(l)}=U_qf_{(l)}$, $g_{(c)}=U_qf_{(c)}$ and $g_{(s)}=U_qf_{(s)}$.
In particular, if $\lambda_1= 2\lambda_{\gamma(f)}$ then $g_{(s)}=0$.

If $f_{(c)}=0$, then $g=g_{(l)}+g_{(s)}$,
thus we have
$$e^{\lambda_1(t+\tau)/2}\Big(\langle g,X_{t+\tau}\rangle-F_{t+\tau}(g_{(s)})\Big)=Y^{1,g_{(l)}}_t(\tau)+Y^{3,g_{(s)}}_t(\tau).$$
Using the convergence of the first, second and fourth components in Theorem \ref{Thm1},
we get for any nonzero $\mu\in {\cal M}_F(E)$, it holds under $\P_{\mu}$ that, as $t\to\infty$,
\begin{equation}\label{large1}
   \left(W_t,~e^{\lambda_1(t+\cdot)/2}\Big(\langle g,X_{t+\cdot}\rangle-F_{t+\cdot}(g_{(s)})\Big)\right)
  \stackrel{d}{\rightarrow}(W_\infty,~G^{1,g_{(l)}}+G^{3,g_{(s)}}),
\end{equation}
where $G^{1,g_{(l)}}+G^{3,g_{(s)}}$ is a continuous Gaussian process with mean 0 and covariance function
\begin{eqnarray*}
&& E\Big[(G^{1,g_{(l)}}(\tau_1)+G^{3,g_{(s)}}(\tau_1))(G^{1,g_{(l)}}(\tau_2)+G^{3,g_{(s)}}(\tau_2))\Big]\\
 &=&\sigma_{g_{(l)},\tau_2-\tau_1}+\eta_{\tau_1,\tau_2}(g_{(l)},g_{(s)})+\beta_{g_{(s)},\tau_2-\tau_1},\quad 0\le \tau_1\le\tau_2.
\end{eqnarray*}

If $f_{(c)}\ne 0$, then
$$t^{-1/2}e^{\lambda_1(t+\tau)/2}\Big(\langle g,X_{t+\tau}\rangle-F_{t+\tau}(g_{(s)})\Big)
=t^{-1/2}\Big(Y^{1,g_{(l)}}_t(\tau)+Y^{3,g_{(s)}}_t(\tau)\Big)+Y^{2,g_{(c)}}_t(\tau).$$
By \eqref{large1}, we get
$$
t^{-1/2}\Big(Y^{1,g_{(l)}}_t(\cdot)+Y^{3,g_{(s)}}_t(\cdot)\Big)
  \stackrel{d}{\rightarrow} 0.
$$
Thus using the convergence of the first and third components in Theorem \ref{Thm1}, we get
$$
\left(W_t,~t^{-1/2}e^{\lambda_1(t+\cdot)/2}\Big(\langle g,X_{t+\cdot}\rangle-F_{t+\cdot}(g_{(s)})\Big) \right)
  \stackrel{d}{\rightarrow}(W_\infty,~G^{2,g_{(c)}}),
$$
where $G^{2,g_{(c)}}\sim \mathcal{N}(0,\rho_{g_{(c)}}^2)$ is a constant process.
Moreover, $W_\infty$ and $G^{2,g_{(c)}}$ are independent.
Note that, if $\lambda_1= 2\lambda_{\gamma(f)}$, then $F_{t+\cdot}(g_{(s)})=0$,
and thus we have $
\left(W_t,~t^{-1/2}e^{\lambda_1(t+\cdot)/2}\langle g,X_{t+\cdot}\rangle\right)
  \stackrel{d}{\rightarrow}(W_\infty,~G^{2,g_{(c)}})
$.
\end{remark}

\section{Preliminaries}
In this section, we give some useful results and facts.
In the remainder of this paper we will use the following notation: for two positive functions $f$ and
$g$ on $E$, $f(x)\lesssim g(x)$ means that there exists a constant $c > 0$ such that $f(x) \le cg(x)$ for all $x\in  E.$

In \cite[(2.25)]{RSZ2}, we  have proved that
\begin{equation}\label{1.19}
  \int_{0}^{t_0}T_s(a_{2t_0})(x)\,ds\lesssim a_{t_0}(x)^{1/2}.
\end{equation}

\subsection{Estimates on the moments of $X$}
In this subsection, we will recall some results about the moments of $\langle f,X_t\rangle$.
The first result is \cite[Lemma 2.1]{RSZ2}.

\begin{lemma}\label{lem:expansion}
 For any $f\in L^2(E,m)$, $x\in E$ and $t>0$, we have
\begin{equation}\label{1.17}
T_t f(x)=\sum_{k=\gamma(f)}^\infty e^{-\lambda_k t}\sum_{j=1}^{n_k}a_j^k\phi^{(k)}_j(x)
\end{equation}
and
\begin{equation}\label{1.25}
  \lim_{t\to\infty}e^{\lambda_{\gamma(f)}t}T_t f(x)=\sum_{j=1}^{n_{\gamma(f)}}a_j^{\gamma(f)}\phi_j^{(\gamma(f))}(x),
\end{equation}
where the series in \eqref{1.17} converges absolutely and uniformly in any compact subset of $E$.
Moreover, for any $t_1>0$,
\begin{eqnarray}
   &&\sup_{t>t_1}e^{\lambda_{\gamma(f)}t}|T_t f(x)|
   \le e^{\lambda_{\gamma(f)}t_1}\|f\|_2\left(\int_E a_{t_1/2}(x)\,m(dx)\right)a_{t_1}(x)^{1/2},\label{1.36}\\
  &&\sup_{t>t_1}e^{(\lambda_{\gamma(f)+1}-\lambda_{\gamma(f)})t}\left|e^{\lambda_{\gamma(f)}t}T_t f(x)-f^*(x)\right|
  \le e^{\lambda_{\gamma(f)+1}t_1}\|f\|_2\left(\int_E a_{t_1/2}(x)\,m(dx)\right)(a_{t_1}(x))^{1/2},\nonumber\\
  &&\label{1.43}
\end{eqnarray}
 where $f^*=\sum_{j=1}^{n_{\gamma(f)}}b^{\gamma(f)}_j\phi_j^{(\gamma(f))}.$
\end{lemma}

We now recall the second moments of the superprocess $\{X_t: t\ge 0\}$ (see, for example, \cite{RSZ3}):
for $f\in L^2(E,m)\cap L^4(E,m)$ and $\mu\in\mathcal{M}_F(E)$, we have for any $t >0$,
\begin{equation}\label{1.9}
  \P_{\mu}\langle f,X_t\rangle^2=\left(\P_{\mu}\langle f,X_t\rangle\right)^2+\int_E\int_0^tT_{s}[A(T_{t-s}f)^2](x)\,ds\mu(dx).
\end{equation}
Thus,
\begin{equation}\label{1.13}
   {\V}{\rm ar}_{\mu}\langle f,X_t\rangle=\langle{\V}{\rm ar}_{\delta_\cdot}\langle f,X_t\rangle, \mu\rangle
   =\int_E\int_0^tT_{s}[A(T_{t-s}f)^2](x)\,ds\mu(dx),
\end{equation}
where $\mathbb{V}{\rm ar}_{\mu}$ stands for the variance under $\P_{\mu}$.
Moreover, for $f\in L^2(E,m)\cap L^4(E,m)$,
\begin{equation}\label{1.4}
  \V{\rm ar}_{\delta_x}\langle f,X_t\rangle\le e^{Kt}T_t(f^2)(x)\in L^2(E,m).
\end{equation}

The next result is \cite[Lemma 2.6]{RSZ3}.

\begin{lemma} \label{lem:2.2}
Assume that $f\in L^2(E,m)\cap L^4(E,m)$.
\begin{description}
  \item{(1)}
   If $\lambda_1<2\lambda_{\gamma(f)}$, then for any $x\in E$,
  \begin{equation}\label{3.2}
   \lim_{t\to \infty}e^{\lambda_1t}{\V}{\rm ar}_{\delta_x}\langle f,X_t\rangle
  =\sigma^2_f \phi_1(x).
  \end{equation}
  Moreover, for $(t, x)\in (3t_0, \infty)\times E$, we have
  \begin{equation}\label{2.10}
     e^{\lambda_1t}{\V}{\rm ar}_{\delta_x}\langle f,X_t\rangle\lesssim a_{t_0}(x)^{1/2}.
  \end{equation}

  \item{(2)}
  If $\lambda_1=2\lambda_{\gamma(f)}$, then for any $(t, x)\in (3t_0, \infty)\times E$,
  \begin{equation}\label{1.49}
   \left|t^{-1}e^{\lambda_1 t}{\V}{\rm ar}_{\delta_x}\langle f,X_t\rangle-\rho_{f^*}^2\phi_1(x)\right|
   \lesssim t^{-1}a_{t_0}(x)^{1/2},
  \end{equation}
  where $f^*=\sum_{j=1}^{n_{\gamma(f)}}b^{\gamma(f)}_j\phi_j^{(\gamma(f))}.$

  \item{(3)}
If $\lambda_1>2\lambda_{\gamma(f)}$, then for any $x\in E$,
  \begin{equation}\label{2.23}
       \lim_{t\to \infty}e^{2\lambda_{\gamma(f)}t}{\V}{\rm ar}_{\delta_x}\langle f,X_t\rangle=\int_0^\infty e^{2\lambda_{\gamma(f)}s}T_s(A(f^*)^2)(x)\,ds.
  \end{equation}
   Moreover,
  for any $(t, x)\in (3t_0, \infty)\times E$,
  \begin{equation}\label{1.60}
   e^{2\lambda_{\gamma(f)}t}\P_{\delta_x}\langle f,X_t\rangle^2\lesssim a_{t_0}(x)^{1/2}.
  \end{equation}
\end{description}
\end{lemma}

\subsection{Excursion measures of $X$}

We use $\mathbb{D}$ to denote the space of $\mathcal{M}_F({E})$-valued
right continuous functions $t\mapsto \omega_t$ on $(0, \infty)$ having zero as a trap.
We use $(\mathcal{A},\mathcal{A}_t)$ to denote the natural $\sigma$-algebras on
$\mathbb{D}$ generated by the coordinate process.

It is known (see \cite[Section 8.4]{Li11}) that
one can associate with $\{\P_{\delta_x}:x\in E\}$ a family
of $\sigma$-finite measures $\{\N_x:x\in E\}$ defined on $(\D,\mathcal{A})$ such  that $\mathbb{N}_x(\{0\})=0$,
\begin{equation}
 \int_{\D}(1- e^{-\langle f, \omega_{t} \rangle})\mathbb{N}_x(d\omega)
= -\log \mathbb{P}_{\delta_x}(e^{-\langle f, X_{t} \rangle}) ,
\quad f\in {\cal B}^+_b(E),\ t> 0,
\label{DK}
\end{equation}
and, for every $0<t_1<\cdots<t_n<\infty$, and nonzero $\mu_1,\cdots,\mu_n\in {\cal M}_F(E)$,
\begin{equation}\label{TN}
  \mathbb{N}_x(\omega_{t_1}\in d\mu_1,\cdots,\omega_{t_n}\in d\mu_n)
 = \mathbb{N}_x(\omega_{t_1}\in d\mu_1)\P_{\mu_1}(X_{t_2-t_1}\in d\mu_2)\cdots \P_{\mu_{n-1}}(X_{t_n-t_{n-1}}\in d\mu_n).
\end{equation}
For earlier work on excursion measures of superprocesses, see \cite{elk-roe, Li03, E.B2.}.

For any $\mu\in {\cal M}_F(E)$, let $N(d\omega)$
be a Poisson random measure on the space $\mathbb{D}$
with intensity $\int_ E\mathbb{N}_x(d\omega)\mu(dx)$, in a probability space $(\widetilde\Omega, \widetilde{\cal F}, \mathbf{P}_\mu)$.
We define another process $\{\Lambda_t: t\ge 0\}$ by $\Lambda_0=\mu$ and
$$\Lambda_t:=\int_{\mathbb{D}}\omega_tN(d\omega),\quad t>0.$$
Let $\widetilde{\mathcal{F}}_t$ be the $\sigma$-algebra generated by the random variables $\{N(A):A\in \mathcal{A}_t\}$.
Then, $\{\Lambda,(\widetilde{\mathcal{F}}_t)_{t\ge 0}, \mathbf{P}_{\mu}\}$ has the same law as
$\{X,(\mathcal{G}_t)_{t\ge 0},\P_{\mu}\}$, see \cite[Theorem 8.24]{Li11}.
Thus,
\begin{equation}\label{cf}
  \mathbb{P}_\mu\left[\exp\left\{i\theta\langle f, X_{t+s}\rangle\right\}|X_t\right]
  =\mathbf{P}_{X_t}\left[\exp\left(i\theta\langle f,\Lambda_s^t\rangle\right)\right]
  =\exp\left\{\int_{E}\int_{\mathbb{D}}(e^{i\theta\langle f, \omega_s\rangle}-1)\mathbb{N}_x(d\omega)X_t(dx)\right\}.
\end{equation}
The proposition below contains some useful properties of $\mathbb{N}_x$.
The proofs are similar to those in \cite[Corollary 1.2, Proposition 1.1]{E.B2.}.
\begin{prop}\label{Moments of N}
If $\mathbb{P}_{\delta_x}|\langle f, X_t\rangle|<\infty,$ then
\begin{equation}\label{N1}
  \int_{\mathbb{D}}\langle f, \omega_t\rangle\,\mathbb{N}_x(d\omega) =\mathbb{P}_{\delta_x}\langle f, X_t\rangle.
\end{equation}
If $\mathbb{P}_{\delta_x}\langle f, X_t\rangle^2<\infty,$ then
\begin{equation}\label{N2}
  \int_{\mathbb{D}}\langle f, \omega_t\rangle^2\,\mathbb{N}_x(d\omega) =\mathbb{V}ar_{\delta_x}\langle f, X_t\rangle.
\end{equation}
\end{prop}

\bigskip

\subsection{Potential functions}\label{subs:cad}

Recall that $q>\max\{K, -2\lambda_1\}$.
For any $x\in E$ such that $U_q|f|(x)<\infty$, we have
\begin{equation}\label{U_q}
  U_qf(x)=\int_0^\infty e^{-qs}T_sf(x)\,ds.
\end{equation}

\begin{lemma}
If $f\in L^2(E,m)$,
then for any $\mu\in {\cal M}_F(E)$,
\begin{equation}\label{cont2}
  \P_{\mu}\left\{ \langle U_q|f|,X_t\rangle<\infty, \forall\, t\ge 0\right\}=1.
\end{equation}
Moreover,
$\langle U_qf,X_t\rangle$ is finite and right continuous, $\P_{\mu}$-a.s.
\end{lemma}
\textbf{Proof:}
First, we claim that, if  $f$ is nonnegative and bounded, $e^{-qt}\langle U_qf,X_t\rangle$ is a
nonnegative right continuous supermartingale with respect to $\{\mathcal{G}_t:t\ge 0\}$.
In fact, since
$T_tf(x)\le \|f\|_\infty e^{Kt}$, we have
$$U_qf(x)\le \|f\|_\infty \int_{0}^\infty e^{-qt}e^{Kt}\,dt=(q-K)^{-1}\|f\|_\infty<\infty.$$
Since $T_tf(x)$ is continuous, by the dominated convergence theorem, we get that
$U_qf$ is continuous. Thus, $U_qf$ is a bounded and continuous function on $E$.
Since $X$ is a right continuous process in ${\cal M}_F(E)$, we get that
$t\mapsto\langle U_qf,X_t\rangle$ is right continuous.
By Fubini's theorem, we have, for any $x\in E$ and $t\ge 0$,
$$T_t[U_qf](x)=\int_0^\infty e^{-qs}T_{t+s}f(x)\,ds= e^{qt} \int_t^\infty e^{-qs}T_{s}f(x)\,ds\le e^{qt}U_qf(x).$$
By the Markov property of $X$, we have, for $t>s$,
$$\P_{\mu}(e^{-qt}\langle U_qf,X_t\rangle |\mathcal{G}_s)=e^{-qt}\langle T_{t-s}(U_qf),X_s\rangle\le e^{-qs}\langle U_qf,X_s\rangle.$$
Thus, $e^{-qt}\langle U_qf,X_t\rangle$ is a supermartingale.

Now, if  $f\in L^2(E,m)$ is nonnegative, then $f_M(x):=f(x)\textbf{1}_{f\le M}(x)$ is bounded.
So $e^{-qt}\langle U_q(f_M),X_t\rangle$ is a nonnegative right continuous
supermartingale with respect to $\{\mathcal{G}_t:t\ge 0\}$,
and, as $M\to \infty$.
$$\forall\, t\ge 0: e^{-qt}\langle U_q(f_M),X_t\rangle \uparrow e^{-qt}\langle U_qf,X_t\rangle.$$
Since $U_qf\in L^2(E,m)$, $\P_{\mu}\langle U_qf,X_t\rangle=\langle T_t(U_qf),\mu\rangle<\infty.$
Thus, by \cite[Section 1.4, Theorem 5]{Chung}, $e^{-qt}\langle U_qf,X_t\rangle$ is a right continuous supermartingale.
By \cite[Section 1.4, Corollary 1]{Chung}, $e^{-qt}\langle U_qf,X_t\rangle$ is bounded on each finite interval,
$\P_{\mu}$-a.s., which implies that
for any $N>0$,
$$\P_{\mu}\left(e^{-qt}\langle U_qf,X_t\rangle<\infty, \quad t\in[0,N]\right)=1.$$
Thus, we have
$$\P_{\mu}\left(\langle U_qf,X_t\rangle<\infty, \quad t\in[0,\infty)\right)=1.$$

Finally, we consider  general $f\in L^2(E,m)$. Let
$$\Omega_0:=\left\{\langle U_q|f|,X_t\rangle<\infty, \forall t\ge 0\right\}
\cap \left\{\omega:\langle U_q(f^+),X_t(\omega)\rangle\mbox{ and }\langle U_q(f^-),X_t(\omega)\rangle\mbox{ are right continuous}\right\}.$$
We have proved that, for any $\mu\in {\cal M}_F(E)$,
$\P_{\mu}\left(\Omega_0\right)=1.$
It follows that, for $\omega\in \Omega_0$,
$$\langle U_qf,X_t(\omega)\rangle=\langle U_q(f^+),X_t(\omega)\rangle-\langle U_q(f^-),X_t(\omega)\rangle$$
is well defined and right continuous.
The proof is now complete. \hfill$\Box$

\subsection{Martingale problem of $X$}
In this subsection, we recall the martingale problem of superprocesses.
For more details, see, for instance, \cite[Chapter 7]{Li11}.

For our superprocess $X$,
there exists a worthy
$(\mathcal{G}_t)$-martingale measure $\{M_t(B): t\ge 0; B\in\mathcal{B}(E)\}$
with covariation measure
\begin{equation}\label{5.2}
  \nu(ds,dx,dy):=ds\int_EA(z)\delta_z(dx)\delta_z(dy)X_s(dz)
\end{equation}
such that for $t\ge 0$, $f\in \mathcal{B}_b(E)$ and $\mu\in \mathcal{M}_F(E)$, we have, $\P_{\mu}$-a.s.,
\begin{equation}\label{5.1}
  \langle f,X_t\rangle=\langle T_tf,\mu\rangle+\int_0^t\int_E T_{t-s}f(z)\,M(ds,dz).
\end{equation}
Let $\mathcal{L}_\nu^2(E)$ be the space of two-parameter predictable processes $h_s(x)$
such that for all $T>0$ and $\mu\in \mathcal{M}_F(E)$,
\begin{eqnarray*}
 &&\P_{\mu}\left[\int_0^T\int_{E^2}h_s(x)h_s(y)\,\nu(ds,dx,dy)\right]
 =\P_{\mu}\left[\int_0^T\int_E A(z)h_s(z)^2 X_s(dz)\,ds\right]\\
 &=&\int_E\int_0^TT_s [Ah_s^2 ](z) \,ds\mu(dz)<\infty.
\end{eqnarray*}
Then, for $h\in \mathcal{L}_\nu^2(E)$,
$$M_t(h):=\int_0^t\int_E h_s(z) M(ds,dz)$$
is well defined and it is a square-integrable cadlag
$\mathcal{G}_t$-martingale under $\P_\mu$, for each $\mu\in \mathcal{M}_F(E)$, with
\begin{equation}\label{5.3}
  \langle M(h)\rangle_t=\int_0^t\langle Ah_s^2,X_s\rangle \,ds.
\end{equation}
For $f\in L^2(E,m)\cap L^4(E,m)$ and $\mu\in \mathcal{M}_F(E)$, we have
$$\int_E\int_0^tT_s [A(T_{t-s}f)^2 ](z) \,ds\mu(dz)=\V ar_{\mu}\langle f,X_t\rangle<\infty,$$
which implies that
$$\int_0^t\int_E T_{t-s}f(z)\,M(ds,dz)$$
is well defined.
Now, using a routine limit argument, we can show that \eqref{5.1} holds for
all $f\in L^2(E,m)\cap L^4(E,m)$ and $\mu\in \mathcal{M}_F(E)$.

For $f\in L^2(E,m)\cap L^4(E,m)$, $U_qf\in L^2(E,m)\cap L^4(E,m)$.
By \eqref{5.1}, for $t>0$ and $\mu\in \mathcal{M}_F(E)$, we have, $\P_{\mu}$-a.s.,
\begin{eqnarray}
  \langle U_qf,X_{t}\rangle
  &=&\langle T_{t}(U_qf),\mu\rangle+\int_0^{t}\int_E T_{t-s}(U_qf)(z)M(ds,dz)\nonumber\\
  &=& \langle T_{t}(U_qf),\mu\rangle
  +\int_0^{t}\int_E \int_{0}^\infty e^{-qu}T_{u+t-s}f(z)\,duM(ds,dz)\nonumber\\
  &=& \langle T_{t}(U_qf),\mu\rangle
  +e^{qt}\int_0^{t}\int_E \int_{t}^\infty e^{-qu}T_{u-s}f(z)\,duM(ds,dz)\nonumber\\
  &=&\langle T_{t}(U_qf),\mu\rangle
  +e^{qt} \int_{t}^\infty e^{-qu}\,du\int_0^{t}\int_ET_{u-s}f(z)M(ds,dz)\nonumber\\
  &:=&J_1^f(t)+e^{qt}J_2^f(t),\label{5.6}
\end{eqnarray}
where the fourth equality follows from the stochastic Fubini's theorem for martingale measures
(see, for instance, \cite[Theorem 7.24]{Li11}).
Thus,  for $t>0$ and $\mu\in \mathcal{M}_F(E)$,
\begin{equation}\label{4.3}
  \P_{\mu}\left(\langle U_qf,X_t\rangle=J_1^f(t)+e^{qt}J^f_2(t)\right)=1.
\end{equation}
For any $u>0$ and $0\le T\le u$, we define
$$M^{(u)}_T:=\int_{0}^{T}\int_E T_{u-s}f(x)M(ds,dx).$$
Then, for any $\mu\in \mathcal{M}_F(E)$, $\{M^{(u)}_T,0\le T\le u\}$ is a
cadlag square-integrable martingale under $\P_\mu$ with
\begin{equation}\label{5.22}
  \langle M^{u}\rangle_T=\int_0^T \langle A(T_{u-s}f)^2, X_s\rangle\,ds.
\end{equation}
Note that
\begin{equation}\label{6.10}
  \P_{\mu}(M^{(u)}_u)^2=\P_{\mu}\langle M^{u}\rangle_u=\V ar_{\mu}\langle f,X_u\rangle.
\end{equation}

\begin{lemma}\label{lem:cad}
If $f\in L^2(E,m)\cap L^4(E,m)$ and $\mu\in \mathcal{M}_F(E)$, then
$t\mapsto\langle U_qf,X_t\rangle$ is a cadlag process on $[0,\infty)$, $\P_\mu$-a.s.
Moreover,
\begin{equation}\label{version}
  \P_{\mu}\left(\langle U_qf,X_t\rangle=J_1^f(t)+e^{qt}J_2^f(t), \forall t>0\right)=1.
\end{equation}
\end{lemma}
\textbf{Proof:}
Since $\langle U_qf,X_t\rangle$ is right continuous, $\P_{\mu}$-a.s.,
in light of \eqref{4.3}, to prove \eqref{version},
it suffices to prove that $J_1^f(t)$ and $J_{2}^f(t)$ are all cadlag in $(0,\infty)$, $\P_{\mu}$-a.s..

For $J_1^f(t)$, by Fubini's theorem, for $t>0$,
$$J_1^f(t)=e^{qt}\int_{t}^\infty e^{-qs}\langle T_sf,\mu\rangle\,ds.$$
Thus, it is easy to see that $J_1^f(t)$ is continuous in $t\in (0,\infty)$.

Now, we consider $J_{2}^f(t)$.
We claim that, for any $t_1>0$,
\begin{equation}\label{4.27'}
  \P_{\mu}\left( J_{2}^f(t) \mbox{ is cadlag in } [t_1,\infty)\right)=1.
\end{equation}
By the definition of $J_{2}^f$, for $t\ge t_1 $,
\begin{equation}\label{4.7}
  J_{2}^f(t)=\int_{t_1}^{\infty}e^{-qu}M^{(u)}_{t}\textbf{1}_{t<u}\,du.
\end{equation}
Since $t\mapsto M^{(u)}_{t}\textbf{1}_{t<u}$ is right continuous,
by the dominated convergence theorem, to prove \eqref{4.27'}, it suffices to show that
\begin{equation}\label{4.8}
  \P_{\mu}\left(\int_{t_1}^{\infty}e^{-qu}\sup_{t\ge t_1}\left(|M^{(u)}_{t}|\textbf{1}_{t<u}\right)\,du<\infty\right)=1.
\end{equation}
By  the $L_p$ maximum inequality and \eqref{6.10}, we have
\begin{eqnarray}\label{4.10}
  &&\P_{\mu}\left(\int_{t_1}^{\infty}e^{-qu}\sup_{t\ge t_1}\left(|M^{(u)}_{t}|\textbf{1}_{t<u}\right)\,du\right)
  \le 2\int_{t_1}^\infty e^{-qu}\sqrt{\P_{\mu}\left|M^{(u)}_{u}\right|^2}\,du\nonumber\\
  &=&2\int_{t_1}^\infty e^{-qu}\sqrt{\int_E\V ar_{\delta_x}\langle f,X_u\rangle\,\mu(dx)}\,du.
\end{eqnarray}
By \eqref{1.4} and \eqref{1.36}, we have, for $u>t_1$,
$$\int_E\V ar_{\delta_x}\langle f,X_u\rangle\,\mu(dx)\le e^{Ku}\int_ET_u(f^2)(x)\,\mu(dx)
\lesssim  e^{Ku}e^{-\lambda_1u}\int_Ea_{t_1}(x)^{1/2}\,\mu(dx).$$
Since $a_{t_1}(x)$ is continuous in $E$ and $\mu$ has compact support,
it follows that $\int_E a_{t_1}(z)^{1/2}\mu(dz)<\infty.$
Thus, by \eqref{4.10}, we have
\begin{eqnarray*}
  &&\P_{\mu}\left(\int_{t_1}^{\infty}e^{-qu}\sup_{t\ge t_1}\left(|M^{(u)}_{t}|\textbf{1}_{t<u}\right)\,du\right)
  \lesssim  \int_{t_1}^\infty e^{-qu}e^{(K-\lambda_1)u/2}\,du\sqrt{\int_E a_{t_1}(x)^{1/2}\,\mu(dx)}<\infty.
\end{eqnarray*}
Now \eqref{4.8} follows immediately.
Since $t_1>0$ are arbitrary, we have
\begin{equation}\label{4.9}
  \P_{\mu}\left( J_{2}^f(t) \mbox{ is cadlag in } (0,\infty)\right)=1.
\end{equation}
\hfill$\Box$

\section{Proof of the main result}

Suppose that $(X^n)_{n\ge 0}$ and $X$ are all $\D(\R^d)$-valued random variables and  $D$ is a subset of $\R_+$.
If for any $k\ge 1$ and any $t_1, \dots t_k\in D$,
$$
(X_{t_1}^n,X_{t_2}^n,\cdots,X_{t_k}^n)\stackrel{d}{\rightarrow}(X_{t_1},\cdots,X_{t_k}), \quad \mbox{ as } n\to\infty,
$$
then we write
$$
X^n\stackrel{\mathcal{L}(D)}{\longrightarrow}X, \quad \mbox{ as }  n\to\infty.
$$
It is known (see, for example, \cite[Chapter VI, 3.20]{J.J.}) that,
$X^n\stackrel{d}{\longrightarrow}X$ in $\D(\R^d)$ as $n\to \infty$ if and only if
\begin{enumerate}
  \item [(i)] $(X^n)_{n\ge 0}$ is tight in  $\D(\R^d)$,
  \item [(ii)] $X^n\stackrel{\mathcal{L}(D)}{\longrightarrow}X$ as $n\to \infty$ for some dense subset $D$ of $\R_+$.
\end{enumerate}

\subsection{Finite dimensional convergence}
The following lemma is a generalization of \cite[Remark 1.3]{RSZ3}.
\begin{lemma}\label{lem3}
If $f\in L^2(E,m)$ is nonnegative and $\mu\in\mathcal{M}_F(E)$, then
\begin{equation}\label{6.15}
  e^{\lambda_1t}\langle f,X_t\rangle \to (f,\phi_1)_m W_\infty, \mbox{ in }L^1(\P_{\mu}).
\end{equation}
\end{lemma}
\textbf{Proof:}
For any $M>0$, let $f_M(x):=f(x){\bf{1}}_{|f(x)|<M}$ and $\hat{f}_M:=f-f_M$.
It is obvious that $f_M\in L^2(E,m)\cap L^4(E,m)$.
In \cite[Remark 1.3]{RSZ3}, we have proved that
\begin{equation}\label{6.16}
  \lim_{t\to\infty}\P_{\mu}\left|e^{\lambda_1t}\langle f_M,X_t\rangle -
  (f_M,\phi_1)_m W_\infty\right|=0.
\end{equation}
For $t>t_0$, by \eqref{1.36}, we have
\begin{eqnarray}\label{6.17}
  \P_{\mu}\left|e^{\lambda_1t}\langle \hat{f}_M,X_t\rangle - (\hat{f}_M,\phi_1)_m W_\infty\right|
  &\le& e^{\lambda_1t}\P_{\mu}\langle \hat{f}_M,X_t\rangle +
  |(\hat{f}_M,\phi_1)_m |\P_{\mu}(W_\infty)\nonumber\\
  &=& e^{\lambda_1t}\langle T_t\hat{f}_M,\mu\rangle+
  |(\hat{f}_M,\phi_1)_m |\P_{\mu}(W_\infty)\nonumber\\
 & \lesssim & \|\hat{f}_M\|_2.
\end{eqnarray}
By \eqref{6.16} and \eqref{6.17}, we have
\begin{equation}\label{6.18}
  \limsup_{t\to\infty}\P_{\mu}\left|e^{\lambda_1t}\langle f,X_t\rangle -
   (f,\phi_1)_m W_\infty\right|
  \lesssim \|\hat{f}_M\|_2.
\end{equation}
Letting $M\to \infty$, we arrive at \eqref{6.15}.\hfill$\Box$
\bigskip

Recall that
\begin{equation*}
  H_t^{k,j}:=e^{\lambda_k t}\langle\phi_j^{(k)}, X_t\rangle,\quad t\geq 0,
\end{equation*}
and
for $g(x)=\sum_{k: 2\lambda_k<\lambda_1}\sum_{j=1}^{n_k}b_j^k\phi_j^{(k)}(x)$, $x\in E$,
$$F_t(g):=\sum_{k: 2\lambda_k<\lambda_1}\sum_{j=1}^{n_k}e^{-\lambda_kt}b_j^kH_\infty^{k,j},$$
where $H_\infty^{k,j}$ is the martingale limit of $H_t^{k,j}$. And recall that
$$
I_ug(x):=\sum_{k: 2\lambda_k<\lambda_1}\sum_{j=1}^{n_k}e^{\lambda_ku}b_j^k\phi_j^{(k)}(x),\quad x\in E.
$$
It is easy to see that $I_{s+t}g=I_s(I_tg)$ and $T_u(I_ug)=I_u(T_ug)=g$.
Thus, we have, as $u\to\infty$,
\begin{equation}\label{5.24}
  \langle I_ug,X_{t+u}\rangle\to F_t(g),\quad \P_\mu\mbox{-a.s.}
\end{equation}

Define
$$
\widetilde{H}_t^{k,j}(\omega):=e^{\lambda_k t}\langle\phi_j^{(k)}, \omega_t\rangle,\quad t\geq 0,\omega\in\D,$$
and
$$
H_\infty(g)(\omega):=\sum_{k: 2\lambda_k<\lambda_1}\sum_{j=1}^{n_k}b_j^k\widetilde{H}_\infty^{k,j}(\omega).
$$
It follows from \cite[Lemma 3.1]{RSZ3} that the limit
$\widetilde{H}_\infty^{k,j}:=\lim_{t\to \infty}\widetilde{H}_t^{k.j}$
exists $\mathbb{N}_x\mbox{-a.e.}$,
in $L^1(\mathbb{N}_x)$ and in $ L^2(\mathbb{N}_x).$
Then,  as $u\to\infty$,
\begin{equation}\label{5.30}
 \langle I_ug,\omega_u\rangle\to H_\infty(g)(\omega),\quad \mathbb{N}_x\mbox{-a.e.},
\quad \mbox{in }L^1(\mathbb{N}_x) \mbox{ and in } L^2(\mathbb{N}_x).
\end{equation}
Since $\mathbb{N}_x\langle I_ug,\omega_u\rangle=\mathbb{P}_{\delta_x}\langle I_ug,X_u\rangle=g(x)$,
we get that
\begin{equation}\label{L1H}
 \mathbb{N}_x(H_\infty(g))=g(x).
\end{equation}
By \eqref{N2} and \eqref{1.13},
we have
\begin{equation}\label{var:Iu}
 \mathbb{N}_x\langle I_ug,\omega_u\rangle^2=\mathbb{V}{\rm ar}_{\delta_x}\langle I_ug,X_u\rangle
  =\int_0^uT_s\left[ A(I_sg)^2\right](x)\,ds,
\end{equation}
which implies that
\begin{equation}\label{L2H}
 \mathbb{N}_x(H_\infty(g))^2=\int_0^\infty
    T_s\left[ A(I_sg)^2\right](x)\,ds.
\end{equation}

The following simple fact will be used later:
\begin{equation}\label{3.20'}
  \left|e^{ix}-\sum_{m=0}^n\frac{(ix)^m}{m!}\right|\leq \min\left(\frac{|x|^{n+1}}{(n+1)!}, \frac{2|x|^n}{n!}\right).
\end{equation}

\begin{lemma}\label{fin-dim}
Assume that $f\in \C_s$, $h\in\C_c$, $g\in\C_l$ and $\mu\in \mathcal{M}_F(E)$.
Suppose that $Y^{1,f}_t$, $Y^{2,h}_t$, and $Y^{3,g}_t$ are defined as in Theorem \ref{Thm1}.
Then, for any $0\le \tau_1\le \tau_2\cdots\le\tau_k$, under $\P_{\mu}$, as $t\to\infty$,
\begin{eqnarray}\label{e:fin}
  &&\Big(W_t,Y^{1,f}_t(\tau_1),\cdots,Y^{1,f}_t(\tau_k),Y^{2,h}_t(\tau_1),\cdots,Y^{2,h}_t(\tau_k),Y^{3,g}_t(\tau_1),\cdots,Y^{3,g}_t(\tau_k)\Big)\nonumber\\
  &\stackrel{d}{\rightarrow}& \Big(W_\infty, \sqrt{W_\infty}G^{1,f}(\tau_1),\cdots, \sqrt{W_\infty}G^{1,f}(\tau_k),
  \sqrt{W_\infty}G^{2,h},\cdots, \sqrt{W_\infty}G^{2,h}, \nonumber\\
  &&\quad\sqrt{W_\infty}G^{3,g}(\tau_1),\cdots, \sqrt{W_\infty}G^{3,g}(\tau_k)\Big).
\end{eqnarray}
 Here $G^{2,h}\sim \mathcal{N}(0,\rho^2_h)$ is a constant process and
 $\Big(G^{1,f}(\tau_1),\cdots, G^{1,f}(\tau_k), G^{3,g}(\tau_1),\cdots, G^{3,g}(\tau_k)\Big)$ is
 an $\R^{2k}$-valued Gaussian random variable, with mean $0$ and covariance
\begin{equation}\label{3.3'}
  E(G^{1,f}(\tau_j)G^{1,f}(\tau_l))=\sigma_{f,\tau_l-\tau_j},\quad \mbox{for }1\le j\le l\le k,
\end{equation}
\begin{equation}\label{3.4'}
  E(G^{3,g}(\tau_j)G^{3,g}(\tau_l))=\beta_{g,\tau_l-\tau_j},\quad \mbox{for }1\le j\le l\le k,
\end{equation}
and
\begin{equation}\label{3.6'}
 E(G^{3,g}(\tau_j)G^{1,f}(\tau_l))=\left\{\begin{array}{ll}\eta_{\tau_j,\tau_l}(f,g),&\quad\mbox{if }1\le j< l\le k,\\
 0.&\quad\mbox{if }1\le l\le j\le k.\end{array}\right.
\end{equation}
Moreover, $W_\infty$, $G^{2,h}$ and $\Big(G^{1,f}(\tau_1),\cdots, G^{1,f}(\tau_k), G^{3,g}(\tau_1),\cdots, G^{3,g}(\tau_k)\Big)$  are independent.
\end{lemma}
\textbf{Proof:}
We put $\theta_{1,0}=\theta_{2,0}=\theta_{3,0}=0$,
 $\tau_0=0$ and $s_j:=\tau_j-\tau_{j-1}$, $j=1,\cdots,k$.
Let
$\theta, \theta_{l,j}\in \R$, $l=1,2,3$, $j=1,\cdots,k$.
Define, for $l=0,\cdots,k$,
$$\tilde{f}_l(x):=\sum_{j=l}^k\theta_{1,l}e^{\lambda_1(\tau_{j}-\tau_{l})/2}T_{\tau_j-\tau_l}f(x),\quad \widehat{g}_l(x):=\sum_{j=0}^l\theta_{3,j}e^{\lambda_1(\tau_j-\tau_l)/2}I_{\tau_l-\tau_j}g(x),$$
and
$$B_l(x):=\tilde{f}_l(x)+\theta_{3,l}g(x)-\widehat{g}_l(x).$$
For $j=1,\cdots,k$, by \eqref{5.24},
\begin{equation}\label{7.10}
  F_{t+\tau_j}(g)=
  \lim_{u\to\infty}\langle I_{u+\tau_k-\tau_j}g,X_{u+t+\tau_k}\rangle.
\end{equation}
Then, by \eqref{7.10}, we get that
\begin{eqnarray*}
  &&\P_{\mu}\exp\Big\{i\theta W_t+\sum_{j=1}^ki\theta_{1,j}Y_{t}^{1,f}(\tau_j)+\sum_{j=1}^ki\theta_{2,j}Y^{2,h}_{t}(\tau_j)+\sum_{j=1}^ki\theta_{3,j}Y^{3,g}_{t}(\tau_j)\Big\}\\
  &=& \P_{\mu}\exp\Big\{i\theta W_t+\sum_{j=1}^k\Big[i\theta_{1,j}Y_{t}^{1,f}(\tau_j)+i\theta_{2,j}Y^{2,h}_{t}(\tau_j)+i\theta_{3,j}e^{\lambda_1(t+\tau_j)/2}\big(\langle
 g,X_{t+\tau_j}\rangle-F_{t+\tau_j}(g)\big)\Big]\Big\} \\
 &=& \lim_{u\to\infty} \P_{\mu}\exp\Big\{i\theta W_t+\sum_{j=1}^k\Big[i\theta_{1,j}Y_{t}^{1,f}(\tau_j)+i\theta_{2,j}Y^{2,h}_{t}(\tau_j)+i\theta_{3,j}e^{\lambda_1(t+\tau_j)/2}\langle g,X_{t+\tau_j}\rangle\Big]\\
 &&\quad-i\big\langle \sum_{j=1}^k\theta_{3,j}e^{\lambda_1(t+\tau_j)/2}I_{u+\tau_k-\tau_j}g,X_{u+t+\tau_k}\big\rangle \Big\}\\
 &=&\lim_{u\to\infty} \P_{\mu}\exp\Big\{i\theta W_t+\sum_{j=1}^kie^{\lambda_1(t+\tau_j)/2}\langle \theta_{1,j}f+t^{-1/2}\theta_{2,j}h+\theta_{3,j}g,X_{t+\tau_j}\rangle-ie^{\lambda_1(t+\tau_k)/2}\langle I_{u}\widehat{g}_k,X_{u+t+\tau_k}\rangle \Big\}\\
 &=& \lim_{u\to\infty} \P_{\mu}\exp\Big\{i\theta W_t+
 \sum_{j=1}^kie^{\lambda_1(t+\tau_j)/2}
 \langle \theta_{1,j}f+t^{-1/2}\theta_{2,j}h+\theta_{3,j}g,X_{t+\tau_j}\rangle -ie^{\lambda_1(t+\tau_k)/2}\langle \widehat{g}_k,X_{t+\tau_k}\rangle\\
 &&\quad+i\langle J_u^{(k)}(t,\cdot),X_{t+\tau_k}\rangle \Big\},
\end{eqnarray*}
where
$$
  J_u^{(k)}(t,x)
 :=
  \int_\D \Big(\exp\big\{-ie^{\lambda_1(t+\tau_k)/2}\langle I_u(\widehat{g}_k),\omega_u \rangle\big\}-1+ie^{\lambda_1(t+\tau_k)/2}\langle I_u(\widehat{g}_k),\omega_u \rangle\Big)\N_x(d\omega).
$$
The last equality above follows from the Markov property of $X$, \eqref{cf} and the fact that
$$\int_{\mathbb{D}}\langle I_u\widehat{g}_k,\omega_u\rangle\mathbb{N}_x(d\omega)=\P_{\delta_x}\langle I_u\widehat{g}_k,X_u\rangle=\widehat{g}_k(x).$$
In the proof of \cite[Theorem 1.4]{RSZ3}, we have proved that
$$\lim_{u\to\infty}\langle J_u^{(k)}(t,\cdot),X_{t+\tau_k}\rangle =\langle J^{(k)}(t,\cdot),X_{t+\tau_k}\rangle,\quad \P_\mu\mbox{-a.s.}$$
where
$$J^{(k)}(t,x):=\int_\D \Big(\exp\big\{-ie^{\lambda_1(t+\tau_k)/2}H_\infty(\widehat{g}_k)\big\}-1+ie^{\lambda_1(t+\tau_k)/2}
H_\infty(\widehat{g}_k)\Big)\N_x(d\omega).$$
Thus, by the dominated convergence theorem, we get that
\begin{eqnarray*}
  &&\P_{\mu}\exp\Big\{i\theta W_t+\sum_{j=1}^ki\theta_{1,j}Y_{t}^{1,f}(\tau_j)+\sum_{j=1}^ki\theta_{2,j}Y^{2,h}_{t}(\tau_j)+\sum_{j=1}^ki\theta_{3,j}Y^{3,g}_{t}(\tau_j)\Big\}\\
  &=&\P_{\mu}\exp\Big\{i\theta W_t+\sum_{j=1}^kie^{\lambda_1(t+\tau_j)/2}
  \langle \theta_{1,j}f+t^{-1/2}\theta_{2,j}h+\theta_{3,j}g,X_{t+\tau_j}\rangle-ie^{\lambda_1(t+\tau_k)/2}\langle \widehat{g}_k,X_{t+\tau_k}\rangle\\
  &&\quad+i\langle J^{(k)}(t,\cdot),X_{t+\tau_k}\rangle \Big\}.
\end{eqnarray*}
It is known  (see \cite[3.44]{RSZ3}) that
$$
\lim_{t\to\infty}\langle J^{(k)}(t,\cdot),X_{t+\tau_k}\rangle
=\exp\Big\{-\frac{1}{2}( \mathbb{N}_{\cdot}(H_\infty(\widehat{g}_k))^2, \phi_1)_m W_\infty\Big\}
\quad \mbox{in $\P_\mu$-probability}.
$$
Since
$B_k(x):=\theta_{1,k}f(x)+\theta_{3,k}g(x) -\widehat{g}_k(x)$,
we have, as $t\to\infty$,
\begin{eqnarray}\label{5.35}
  &&\Big|\P_{\mu}\exp\Big\{i\theta W_t+\sum_{j=1}^ki\theta_{1,j}Y_{t}^{1,f}(\tau_j)+\sum_{j=1}^ki\theta_{2,j}Y^{2,h}_{t}(\tau_j)+\sum_{j=1}^ki\theta_{3,j}Y^{3,g}_{t}(\tau_j)\Big\}\nonumber\\
  &&-\P_{\mu}\exp\Big\{(i\theta-\frac{1}{2}(\mathbb{N}_\cdot(H_\infty(\widehat{g}_k))^2, \phi_1)_m ) W_t
  +\sum_{j=1}^{k-1}ie^{\lambda_1(t+\tau_j)/2}\langle \theta_{1,j}f_j+t^{-1/2}\theta_{2,j}h+\theta_{3,j}g,X_{t+\tau_j}\rangle\nonumber\\
  &&+ie^{\lambda_1(t+\tau_k)/2}
  \langle B_k+t^{-1/2}\theta_{2,k}h,X_{t+\tau_k}\rangle\Big\}\Big|\nonumber\\
  &&\to 0.
\end{eqnarray}
By the Markov property of $X$, we have
\begin{eqnarray*}
 && \P_{\mu}\left[\exp\left\{ie^{\lambda_1(t+\tau_k)/2}\langle B_k+t^{-1/2}\theta_{2,k}h,X_{t+\tau_k}\rangle\right\}|\mathcal{F}_{t+\tau_{k-1}}\right]\\
  &=&\exp\left\{ \left\langle \int_\D \left(\exp\left\{ie^{\lambda_1(t+\tau_k)/2}\langle B_k+t^{-1/2}\theta_{2,k}h,\omega_{s_k}\rangle\right\} -1\right) \N_{\cdot}(d\omega),X_{t+\tau_{k-1}}\right\rangle\right\}\\
  &=&\exp\Big\{ie^{\lambda_1(t+\tau_k)/2}\langle \N_{\cdot}\langle B_k+t^{-1/2}\theta_{2,k}h,\omega_{s_k}\rangle,X_{t+\tau_{k-1}}\Big\}\\
 &&\quad\times\exp\Big\{-\frac{1}{2}e^{\lambda_1(t+\tau_k)}\langle\N_{\cdot}
 \langle B_k,\omega_{s_k}\rangle^2,X_{t+\tau_{k-1}}\rangle\Big\}
  \times\exp\left\{\langle R(t,\cdot),X_{t+\tau_{k-1}}\rangle\right\}\\
  &=:&(I)\times(II)\times(III),
\end{eqnarray*}
where
\begin{eqnarray*}
 R(t,x)&:=&\int_\D \left(\exp\left\{ie^{\lambda_1(t+\tau_k)/2}\langle t^{-1/2}\theta_{2,k}h+B_k,\omega_{s_k}\rangle\right\} -1\right.\\
 &&\left.-ie^{\lambda_1(t+\tau_k)/2}\langle t^{-1/2}\theta_{2,k}h+B_k,\omega_{s_k}\rangle+\frac{1}{2}e^{\lambda_1(t+\tau_k)}\langle B_k,\omega_{s_k}\rangle^2\right)
\N_x(d\omega),\quad x\in E.
\end{eqnarray*}
For part $(I)$,
by the definition of $\widehat{g}_k$, we get that
\begin{equation}\label{hatg}
    \theta_{3,k}g(x)-\widehat{g}_k(x)
  =-\sum_{j=0}^{k-1}\theta_{3,j}e^{\lambda_1(\tau_j-\tau_k)/2}I_{\tau_k-\tau_j}g(x)
  =-e^{-\lambda_1(\tau_{k}-\tau_{k-1})/2}I_{\tau_k-\tau_{k-1}}\widehat{g}_{k-1}(x),\quad x\in E.
\end{equation}
Since $h\in\C_c$, we have $T_{s}h(x)=e^{-\lambda_1s/2}h(x)$.
Thus, for $ x\in E$,
\begin{eqnarray*}
  &&\N_x(\langle B_k+t^{-1/2}\theta_{2,k}h,\omega_{s_k}\rangle)
  =T_{s_k}(B_k+t^{-1/2}\theta_{2,k}h)(x)\nonumber\\
 &=&\theta_{1,k}T_{s_k}f(x)+t^{-1/2}\theta_{2,k}e^{-\lambda_1s_k/2}h(x)-e^{-\lambda_1s_k/2}\widehat{g}_{k-1}(x).
  \end{eqnarray*}
Hence, we have
\begin{equation}\label{5.26}
  (I)=\exp\Big\{ie^{\lambda_1(t+\tau_{k-1})/2}\Big\langle\theta_{1,k}e^{\lambda_1s_k/2}T_{s_k}f+t^{-1/2}\theta_{2,k}h
  -\widehat{g}_{k-1},X_{t+\tau_{k-1}}\Big\rangle\Big\}.
\end{equation}
For part $(II)$, we define for $j=1,\cdots,k$,
\begin{equation}\label{def-Cj}
C_j:=e^{\lambda_1s_j}
(\N_{\cdot}\langle B_j,\omega_{s_j}\rangle^2,\phi_1)_m=e^{\lambda_1s_j}
(\V ar_{\delta_\cdot}\langle B_j,\omega_{s_j}\rangle,\phi_1)_m.
\end{equation}
By Lemma \ref{lem3}, we get that, as $t\to\infty$,
\begin{equation*}
  e^{\lambda_1(t+\tau_k)}\langle\N_{\cdot}\langle B_k,\omega_{s_k}\rangle^2,X_{t+\tau_{k-1}}\rangle\to C_kW_\infty
\end{equation*}
in $\P_\mu$-probability.
Thus,  we get that, as $t\to\infty$,
\begin{equation}\label{5.61}
  (II)\to \exp\large\{ -\frac{1}{2}C_kW_\infty\large\}, \mbox{ in $\P_\mu$-probability.}
\end{equation}
Now, we deal with part $(III)$.
For $x_1,x_2\in\R$, by \eqref{3.20'}, we have
\begin{eqnarray}\label{5.62}
  &&\left|e^{i(x_1+x_2)}-1-i(x_1+x_2)+\frac{1}{2}(x_1)^2\right|\nonumber\\
 &\le&\left|e^{ix_1}-1-ix_1+\frac{1}{2}(x_1)^2\right|+\left|e^{ix_2}-1-ix_2\right|+|e^{ix_1}-1||e^{ix_2}-1|\nonumber\\
 &\le& |x_1|^2(1\wedge\frac{|x_1|}{6})+\frac{1}{2}|x_2|^2+|x_1x_2|.
\end{eqnarray}
Using \eqref{5.62} with $x_1=e^{\lambda_1(t+\tau_k)/2}\langle B_k,\omega_{s_k}\rangle$ and $x_2=\theta_{2,k}t^{-1/2}e^{\lambda_1(t+\tau_k)/2}\langle h,\omega_{s_k}\rangle$, we get
\begin{eqnarray*}
  |R(t,x)|&\le& e^{\lambda_1(t+\tau_k)}\N_x\left[\langle B_k,\omega_{s_k}\rangle^2\left(1\bigwedge\frac{e^{\lambda_1(t+\tau_k)/2}|\langle B_k,\omega_{s_k}\rangle|}{6}\right)\right]\\
  &&\quad+\frac{(\theta_{2,k})^2}{2}t^{-1}e^{\lambda_1(t+\tau_k)}\N_x\langle h,\omega_{s_k}\rangle^2
  +|\theta_{2,k}|t^{-1/2}e^{\lambda_1(t+\tau_k)}\N_x\left|\langle h,\omega_{s_k}\rangle\langle B_k,\omega_{s_k}\rangle\right|\\
  &=&e^{\lambda_1(t+\tau_k)}\Big(\N_x\left[\langle B_k,\omega_{s_k}    \rangle^2\left(1\bigwedge\frac{e^{\lambda_1(t+\tau_k)/2}|\langle B_k,\omega_{s_k}\rangle|}{6}\right)\right]\\
     &&\quad+\frac{(\theta_{2,k})^2}{2}t^{-1}\N_x\langle h,\omega_{s_k}\rangle^2+|\theta_{2,k}|t^{-1/2}\N_x
     \left|\langle h,\omega_{s_k}\rangle\langle B_k,\omega_{s_k}\rangle\right|\Big)\\
  &=:&e^{\lambda_1(t+\tau_k)}U(t,x).
\end{eqnarray*}
Notice that $U(\cdot,x)\downarrow0$, as $t\to\infty$.
Thus, for $t>u$,
$$\limsup_{t\to\infty}e^{\lambda_1(t+\tau_k)}\P_{\mu}\langle U(t,\cdot), X_{t+\tau_{k-1}}\rangle\le \limsup_{t\to\infty}e^{\lambda_1(t+\tau_k)}\langle T_{t+\tau_{k-1}} U(u,\cdot),\mu\rangle
=e^{\lambda_1s_k}(U(u,\cdot),\phi_1)_m\langle\phi_1,\mu\rangle,
$$
where the last equality follows from \eqref{1.25}.
Letting $u\to\infty$, we get that
$$\lim_{t\to\infty}e^{\lambda_1(t+\tau_k)}\P_{\mu}\langle U(t,\cdot), X_{t+\tau_{k-1}}\rangle=0,$$
which implies that
\begin{equation}\label{5.28}
  \lim_{t\to\infty}\P_{\mu}|\langle R(t,\cdot), X_{t+\tau_{k-1}}\rangle|=0.
\end{equation}
Thus, by \eqref{5.26}, \eqref{5.61} and \eqref{5.28}, we have that, as $t\to\infty$,
\begin{eqnarray*}
  &&\left|\P_{\mu}\left[\exp\big\{ie^{\lambda_1(t+\tau_k)/2}\langle B_k+t^{-1/2}\theta_{2,k}h,X_{t+\tau_k}\rangle\big\}|\mathcal{F}_{t+\tau_{k-1}}\right]\right.\\
  &&\left.-\exp\Big\{-\frac{1}{2}C_kW_t+ie^{\lambda_1(t+\tau_{k-1})/2}\Big\langle\theta_{1,k}e^{\lambda_1(s_k)/2}T_{s_k}f
  +t^{-1/2}\theta_{2,k}h
  -\widehat{g}_{k-1},X_{t+\tau_{k-1}}\Big\rangle\Big\}\right|\\
  &&\to0 \qquad \mbox{in $\P_\mu$-probability.}
\end{eqnarray*}
Hence, using the Markov property and the dominated convergence theorem, we get that, as $t\to\infty$,
\begin{eqnarray*}
  &&\left|\P_{\mu}\exp\Big\{(i\theta-\frac{1}{2}(\mathbb{N}_\cdot(H_\infty(\widehat{g}_k))^2, \phi_1)_m ) W_t
  +\sum_{j=1}^{k-1}ie^{\lambda_1(t+\tau_j)/2}\langle \theta_{1,j}f+t^{-1/2}\theta_{2,j}h+\theta_{3,j}g,X_{t+\tau_j}\rangle\right.\nonumber\\
  &&\left.+ie^{\lambda_1(t+\tau_k)/2}
  \langle B_k+t^{-1/2}\theta_{2,k}h,X_{t+\tau_k}\rangle\Big\}-\right.\nonumber\\
  &&\left.\P_{\mu}\exp\Big\{(i\theta-\frac{1}{2}(\mathbb{N}_\cdot(H_\infty(\widehat{g}_k))^2, \phi_1)_m -\frac{1}{2}C_k) W_t
  +\sum_{j=1}^{k-2}ie^{\lambda_1(t+\tau_j)/2}\langle \theta_{1,j}f+t^{-1/2}\theta_{2,j}h+\theta_{3,j}g,X_{t+\tau_j}\rangle\right.\\
  &&\left.+ie^{\lambda_1(t+\tau_{k-1})/2}
  \langle B_{k-1}
  +t^{-1/2}(\theta_{2,k-1}+\theta_{2,k})h,X_{t+\tau_{k-1}}\rangle\Big\}\right|\nonumber\\
  &&\to 0.
\end{eqnarray*}
Repeating the above procedure $k$ times, we obtain that, as $t\to\infty$,
\begin{eqnarray}\label{5.29}
 && \left|\P_{\mu}\exp\Big\{(i\theta-\frac{1}{2}(\mathbb{N}_\cdot(H_\infty(\widehat{g}_k))^2, \phi_1)_m ) W_t
  +\sum_{j=1}^{k-1}ie^{\lambda_1(t+\tau_j)/2}\langle \theta_{1,j}f+t^{-1/2}\theta_{2,j}h+\theta_{3,j}g,X_{t+\tau_j}\rangle\right.\nonumber\\
  &&\left.+ie^{\lambda_1(t+\tau_k)/2}
  \langle B_k+t^{-1/2}\theta_{2,k}h,X_{t+\tau_k}\rangle\Big\}-\right.\nonumber\\
  &&\left.-\P_{\mu}\exp\Big\{(i\theta-\frac{1}{2}(\mathbb{N}_\cdot(H_\infty(\widehat{g}_k))^2,\phi_1)_m-\frac{1}{2}\sum_{j=1}^kC_j)W_t+ie^{\lambda_1t/2}\langle \tilde{f}_0,X_t\rangle+it^{-1/2}e^{\lambda_1t/2}
  \langle\sum_{j=1}^k\theta_{2,j}h,X_t\rangle\Big\}\right|\nonumber\\
 &&\to 0.
\end{eqnarray}
By \cite[Lemma 3.5]{RSZ3}, we have
\begin{eqnarray}\label{5.64}
  &&\lim_{t\to\infty}\exp\Big\{\big(i\theta-\frac{1}{2}(\mathbb{N}_\cdot(H_\infty(\widehat{g}_k))^2,\phi_1)_m-\frac{1}{2}\sum_{j=1}^kC_j\big)W_t+ie^{\lambda_1t/2}\langle \tilde{f}_0,X_t\rangle+it^{-1/2}e^{\lambda_1t/2}
  \langle\sum_{j=1}^k\theta_{2,j}h,X_t\rangle\Big\}\nonumber\\
  &=&\exp\Big\{\Big(i\theta-\frac{1}{2}(\mathbb{N}_\cdot(H_\infty(\widehat{g}_k))^2,\phi_1)_m-\frac{1}{2}\sum_{j=1}^kC_j
  -\frac{1}{2}\sigma_{\tilde{f}_0}^2-\frac{1}{2}\Big(\sum_{j=1}^k\theta_{2,j}\Big)^2\rho_{h}^2\Big)W_\infty\Big\}.
\end{eqnarray}
Thus, by \eqref{5.35}, \eqref{5.29} and \eqref{5.64}, we get
\begin{eqnarray}\label{5.67}
  &&\lim_{t\to\infty}\P_{\mu}\exp\Big\{i\theta W_t+\sum_{j=1}^ki\theta_{1,j}Y_{t}^{1,f}(\tau_j)+\sum_{j=1}^ki\theta_{2,j}Y^{2,h}_{t}(\tau_j)+\sum_{j=1}^ki\theta_{3,j}Y^{3,g}_{t}(\tau_j)\Big\}\nonumber\\
  &=&\exp\Big\{\Big(i\theta-\frac{1}{2}(\mathbb{N}_\cdot(H_\infty(\widehat{g}_k))^2,\phi_1)_m-\frac{1}{2}\sum_{j=1}^kC_j
  -\frac{1}{2}\sigma_{\tilde{f}_0}^2-\frac{1}{2}\Big(\sum_{j=1}^k\theta_{2,j}\Big)\rho_{h}^2\Big)W_\infty\Big\}.
\end{eqnarray}
By the definition of $C_j$ in \eqref{def-Cj}, we have,
\begin{eqnarray*}
  &&(\mathbb{N}_\cdot(H_\infty(\widehat{g}_k))^2,\phi_1)_m+\sum_{j=1}^kC_j
  +\sigma_{\tilde{f}_0}^2\\
  &=&\Big[(\mathbb{N}_\cdot(H_\infty(\widehat{g}_k))^2,\phi_1)_m+\sum_{j=1}^ke^{\lambda_1s_j}(\V ar_{\delta_\cdot}\langle \theta_{3,j}g-\widehat{g}_j,\omega_{s_j}\rangle,\phi_1)_m\Big]\\
  &&+\Big[\sum_{j=1}^ke^{\lambda_1s_j}(\V ar_{\delta_\cdot}\langle \tilde{f}_j,\omega_{s_j}\rangle,\phi_1)_m+\sigma_{\tilde{f}_0}^2\Big]
  +2\sum_{j=1}^ke^{\lambda_1s_j}(\mathbb{C}ov_{\delta_\cdot}(\langle \tilde{f}_j,\omega_{s_j}\rangle,\langle \theta_{3,j}g-\widehat{g}_j,\omega_{s_j}\rangle),\phi_1)_m.
\end{eqnarray*}

In the following, we calculate the three parts separately.
\begin{enumerate}
\item
By \eqref{L2H} and \eqref{hatg}, we  have that, for $j=1,\cdots,k$,
\begin{eqnarray*}
&&( \mathbb{N}_\cdot(H_\infty(\widehat{g}_j))^2,\phi_1)_m=\int_0^\infty e^{-\lambda_1s}( A(I_s\widehat{g}_j)^2,\phi_1)_m\,ds\\
  &=&\int_0^\infty e^{-\lambda_1s}
\left(A(I_s(\theta_{3,j}g+e^{-\lambda_1(\tau_j-\tau_{j-1})/2}I_{\tau_j-\tau_{j-1}}\widehat{g}_{j-1}))^2,\phi_1\right)_m\,ds\\
  &=&\theta_{3,j}^2\beta_{g}^2
  +2\theta_{3,j}\sum_{l=0}^{j-1}\theta_{3,l}\beta_{g,\tau_j-\tau_l}
  +\int_{\tau_j-\tau_{j-1}}^\infty
   e^{-\lambda_1s}(A(I_s\widehat{g}_{j-1})^2,\phi_1)_m\,ds.
\end{eqnarray*}
By \eqref{var:Iu} and \eqref{hatg}, we get that
\begin{eqnarray*}
  &&\V ar_{\delta_\cdot}\langle \theta_{3,j}g-\widehat{g}_j,\omega_{s_j}\rangle=( \V ar_{\delta_\cdot}\langle I_{\tau_j-\tau_{j-1}}\widehat{g}_{j-1},\omega_{\tau_j-\tau_{j-1}}\rangle,\phi_1)_m\nonumber\\
 &=&\int_{0}^{\tau_j-\tau_{j-1}} e^{-\lambda_1s}( A[I_s(\widehat{g}_{j-1})]^2,\phi_1)_m\,ds.
\end{eqnarray*}
Thus, we have, for $j=1,\cdots,k$,
$$
( \mathbb{N}_\cdot(H_\infty(\widehat{g}_j))^2,\phi_1)_m+(\V ar_{\delta_\cdot}\langle \theta_{3,j}g-\widehat{g}_j,\omega_{s_j}\rangle,\phi_1)_m
=\theta_{3,j}^2\beta_{g}^2
+2\theta_{3,j}\sum_{l=0}^{j-1}\theta_{3,l}\beta_{g,\tau_j-\tau_l}
+( \mathbb{N}_\cdot(H_\infty(\widehat{g}_{j-1}))^2,\phi_1)_m.
$$
Summing over $j$ and using the fact that $\widehat{g}_0=0$, we get
\begin{equation}\label{5.66}
  (\mathbb{N}_\cdot(H_\infty(\widehat{g}_k))^2,\phi_1)_m+\sum_{j=1}^k(\V ar_{\delta_\cdot}\langle \theta_{3,j}g-\widehat{g}_j,\omega_{s_j}\rangle,\phi_1)_m=
\sum_{j=1}^k\theta_{3,j}^2\beta_{g}^2
+2\sum_{j=1}^k\sum_{l=0}^{j-1} \theta_{3,j}\theta_{3,l}\beta_{g,\tau_j-\tau_l}.
\end{equation}
\item
Since $\tilde{f}_j=\theta_{1,j}f+e^{\lambda_1(\tau_{j+1}-\tau_j)/2}T_{\tau_{j+1}-\tau_j}\tilde{f}_{j+1}
=\theta_{1,j}f+\sum_{l=j+1}^k\theta_{1,l}e^{\lambda_1(\tau_{l}-\tau_j)/2}T_{\tau_{l}-\tau_j}f$,
we have
\begin{eqnarray*}
 \sigma_{\tilde{f}_{j}}^2&=&\int_0^\infty e^{\lambda_1u}(A[T_u\tilde{f}_{j}]^2,\phi_1)_m\,du \\
  &=& \int_0^\infty e^{\lambda_1u}(A[T_u(\theta_{1,j}f+e^{\lambda_1(\tau_{j+1}-\tau_{j})/2}T_{\tau_{j+1}-\tau_{j}}\tilde{f}_{j+1})]^2,\phi_1)_m\,du\\
  &=&\theta_{1,j}^2\sigma_f^2+2\sum_{l=j+1}^k\theta_{1,j}\theta_{1,l}\sigma_{f,\tau_l-\tau_{j}}+e^{\lambda_1(\tau_{j+1}-\tau_j)}\int_0^\infty e^{\lambda_1u}(A[T_{u+\tau_{j+1}-\tau_j}\tilde{f}_{j+1}]^2,\phi_1)_m\,du\\
  &=&\theta_{1,j}^2\sigma_f^2+2\sum_{l=j+1}^k\theta_{1,j}\theta_{1,l}\sigma_{f,\tau_l-\tau_j}+\int_{\tau_{j+1}-\tau_j}^\infty e^{\lambda_1u}(A[T_{u}\tilde{f}_{j+1}]^2,\phi_1)_m\,du.
\end{eqnarray*}
By \eqref{1.13}, we have
$$e^{\lambda_1s_{j+1}}(\V ar_{\delta_\cdot}\langle \tilde{f}_{j+1},\omega_{s_{j+1}}\rangle,\phi_1)_m=\int_0^{\tau_{j+1}-\tau_j} e^{\lambda_1u}(A[T_{u}\tilde{f}_{j+1}]^2,\phi_1)_m\,du.$$
Thus, we get, for $j=0,\cdots,k-1$,
$$\sigma_{\tilde{f}_{j}}^2+e^{\lambda_1s_{j+1}}(\V ar_{\delta_\cdot}\langle \tilde{f}_{j+1},\omega_{s_{j+1}}\rangle,\phi_1)_m=\theta_{1,j}^2\sigma_f^2+2\sum_{l=j+1}^k\theta_{1,j}\theta_{1,l}\sigma_{f,\tau_l-\tau_j}
+\sigma_{\tilde{f}_{j+1}}^2.$$
Therefore, summing over $j$ on both sides of the above equality,
 we get
\begin{eqnarray}\label{5.65}
  &&\sum_{j=1}^ke^{\lambda_1s_j}(\V ar_{\delta_\cdot}\langle \tilde{f}_j,\omega_{s_j}\rangle,\phi_1)_m+\sigma_{\tilde{f}_0}^2\nonumber\\
  &=&\sum_{j=0}^{k-1}\theta_{1,j}^2\sigma_f^2+2\sum_{j=0}^{k-1}\sum_{l=j+1}^k\theta_{1,j}\theta_{1,l}\sigma_{f,\tau_l-\tau_j}+\sigma_{\tilde{f}_{k}}^2\nonumber\\
  &=&\sum_{j=1}^{k}\theta_{1,j}^2\sigma_f^2+2\sum_{j=1}^{k-1}\sum_{l=j+1}^k\theta_{1,j}\theta_{1,l}\sigma_{f,\tau_l-\tau_j},
\end{eqnarray}
where the last equality follows from the fact that $\theta_{1,0}=0$ and $\tilde{f}_k=\theta_{1,k}f$.
\item
Since $\tilde{f}_j=\sum_{l=j}^k\theta_{1,l}e^{\lambda_1(\tau_l-\tau_j)/2}T_{\tau_l-\tau_j}f$ and $\theta_{3,j}g-\widehat{g}_j=-\sum_{r=0}^{j-1}\theta_{3,r}e^{-\lambda_1(\tau_j-\tau_r)/2}I_{\tau_j-\tau_r}g$,
we have
\begin{eqnarray*}
 &&e^{\lambda_1s_j}(\mathbb{C}ov_{\delta_\cdot}(\langle \tilde{f}_j,\omega_{s_j}\rangle,\langle \theta_{3,j}g-\widehat{g}_j,\omega_{s_j}\rangle),\phi_1)_m\\
  &=& \int_0^{\tau_j-\tau_{j-1}} e^{\lambda_1u}\Big(AT_u(\tilde{f}_j)T_u(\theta_{3,j}g-\widehat{g}_j),\phi_1\Big)_m\,du\\
  &=&-\sum_{l=j}^k\sum_{r=0}^{j-1}\theta_{1,l}\theta_{3,r}
  e^{\lambda_1(\tau_l+\tau_r-2\tau_j)/2}\int_0^{\tau_j-\tau_{j-1}} e^{\lambda_1u}\Big(AT_{u+\tau_l-\tau_j}fI_{\tau_j-\tau_r-u}g,\phi_1\Big)_m\,du\\
 &=&-\sum_{l=j}^k\sum_{r=0}^{j-1}\theta_{1,l}\theta_{3,r}
  e^{\lambda_1(\tau_l+\tau_r)/2}\int_{\tau_{j-1}}^{\tau_j} e^{-\lambda_1u}\Big(AT_{\tau_l-u}fI_{u-\tau_r}g,\phi_1\Big)_m\,du.
\end{eqnarray*}
Thus, we get that
\begin{eqnarray}\label{5.8}
  &&2\sum_{j=1}^ke^{\lambda_1s_j}(\mathbb{C}ov_{\delta_\cdot}(\langle \tilde{f}_j,\omega_{s_j}\rangle,\langle \theta_{3,j}g-\widehat{g}_j,\omega_{s_j}\rangle),\phi_1)_m\nonumber\\
  &=&-2\sum_{l=1}^k\sum_{r=0}^{l-1}\sum_{j=r+1}^l\theta_{1,l}\theta_{3,r}
  e^{\lambda_1(\tau_l+\tau_r)/2}\int_{\tau_{j-1}}^{\tau_j} e^{-\lambda_1u}\Big(AT_{\tau_l-u}fI_{u-\tau_r}g,\phi_1\Big)_m\,du\nonumber\\
  &=&-2\sum_{l=1}^k\sum_{r=0}^{l-1}\theta_{1,l}\theta_{3,r}
  e^{\lambda_1(\tau_l+\tau_r)/2}\int_{\tau_{r}}^{\tau_l} e^{-\lambda_1u}\Big(AT_{\tau_l-u}fI_{u-\tau_r}g,\phi_1\Big)_m\,du.
\end{eqnarray}
\end{enumerate}
Combining \eqref{5.67}--\eqref{5.8}, we get \eqref{e:fin} immediately.

The proof is now complete. \hfill$\Box$

\begin{remark}\label{rek:4}
{\rm By Lemma \ref{fin-dim}, for any $f\in \C_s$ and $g\in\C_l$, there exists a Gaussian process
$\big(G^{1,U_qf},G^{3,g}\big)$ with mean $0$ and covariance function defined as in Theorem \ref{Thm1}.
Furthermore, the next lemma shows that, this Gaussian process has a continuous version.
Thus, the Gaussian process $\big(G^{1,U_qf},G^{3,g}\big)$ defined in Theorem \ref{Thm1} exists.}
\end{remark}

\begin{lemma}\label{contin}
Assume that $f\in \C_s$ and $g\in\C_l$.
If $\big(G^{1,U_qf}(\tau),G^{3,g}(\tau)\big)_{\tau\ge0}$ is a Gaussian process with mean $0$ and
covariance function defined as in Theorem \ref{Thm1},
then, $\big(G^{1,U_qf},G^{3,g}\big)$ has a continuous version.
\end{lemma}
\textbf{Proof:}
By Kolmogorov's continuity criterion, it suffices to show that,
for any $\tau_2>\tau_1\ge 0$,
\begin{equation}\label{moment}
  E|G^{1,U_qf}(\tau_2)-G^{1,U_qf}(\tau_1)|^4+E|G^{3,g}(\tau_2)-G^{3,g}(\tau_1)|^4\le C|\tau_2-\tau_1|^2,
\end{equation}
where $C$ is a constant.

(1)
Since $G^{1,U_qf}(\tau_2)-G^{1,U_qf}(\tau_1)\sim \mathcal{N}(0,\Sigma(\tau_1,\tau_2))$
with $\Sigma(\tau_1,\tau_2)= E|G^{1,U_qf}(\tau_2)-G^{1,U_qf}(\tau_1)|^2$,
we have
\begin{equation}\label{4.15}
  E|G^{1,U_qf}(\tau_2)-G^{1,U_qf}(\tau_1)|^4=\Sigma(\tau_1,\tau_2)^2E(G^4),
\end{equation}
where $G\sim\mathcal{N}(0,1)$.
In the following, we write $U_qf$ as $f^{(q)}$.
By \eqref{3.3'}, we have
\begin{eqnarray*}
 \Sigma(\tau_1,\tau_2)&=&E|G^{1,U_qf}(\tau_2)-G^{1,U_qf}(\tau_1)|^2 \\
  &=&  2\int_0^\infty e^{\lambda_1s}(A(T_sf^{(q)})^2,\phi_1)_m\,ds
  -2e^{\lambda_1(\tau_2-\tau_1)/2}\int_0^\infty e^{\lambda_1s}(A(T_sf^{(q)})(T_{s+\tau_2-\tau_1}f^{(q)}),\phi_1)_m\,ds\\
   &=& 2\int_0^\infty e^{\lambda_1s}
   \left( A(T_sf^{(q)})(T_s(f^{(q)}-e^{\lambda_1(\tau_2-\tau_1)/2}T_{\tau_2-\tau_1}f^{(q)})),\phi_1\right)_m\,ds\\
   &\le&2K\int_0^\infty e^{\lambda_1s}\left\|(T_sf^{(q)})\left(T_s(f^{(q)}-e^{\lambda_1(\tau_2-\tau_1)/2}T_{\tau_2-\tau_1}f^{(q)})\right)\right\|_2\,ds.
\end{eqnarray*}
We rewrite the last integral above as the sum of integrals over $(0, t_0)$ and $(t_0, \infty).$
For $s>t_0$,
\begin{equation}\label{4.20}
  \left\|(T_sf^{(q)})\left(T_s(f^{(q)}-e^{\lambda_1(\tau_2-\tau_1)/2}T_{\tau_2-\tau_1}f^{(q)})\right)\right\|_2
  \lesssim e^{-2\lambda_{\gamma(f)}s}\|a_{t_0}\|_2\|f^{(q)}\|_2\|f^{(q)}-e^{\lambda_1(\tau_2-\tau_1)/2}T_{\tau_2-\tau_1}f^{(q)}\|_2.
\end{equation}
Thus,
\begin{equation}\label{4.21}
  \int_{t_0}^\infty e^{\lambda_1s}\left\|(T_sf^{(q)})\left(T_s(f^{(q)}-e^{\lambda_1(\tau_2-\tau_1)/2}T_{\tau_2-\tau_1}f^{(q)})\right)\right\|_2\, ds
  \lesssim \|f^{(q)}-e^{\lambda_1(\tau_2-\tau_1)/2}T_{\tau_2-\tau_1}f^{(q)}\|_2.
\end{equation}
For $s\le t_0$, since $\|T_s\|_4\le e^{Ks}$, we have
\begin{eqnarray*}
  &&\|(T_sf^{(q)})\left(T_s(f^{(q)}-e^{\lambda_1(\tau_2-\tau_1)/2}T_{\tau_2-\tau_1}f^{(q)})\right)\|_2
  \le \|T_sf^{(q)}\|_4\,\|T_s(f^{(q)}-e^{\lambda_1(\tau_2-\tau_1)/2}T_{\tau_2-\tau_1}f^{(q)})\|_4\\
   &\le&  e^{2Ks}\|f^{(q)}\|_4\,\|f^{(q)}-e^{\lambda_1(\tau_2-\tau_1)/2}T_{\tau_2-\tau_1}f^{(q)}\|_4.
\end{eqnarray*}
Thus,
\begin{equation}\label{4.22}
  \int_0^{t_0} e^{\lambda_1s}\left\|(T_sf^{(q)})\left(T_s(f^{(q)}-e^{\lambda_1(\tau_2-\tau_1)/2}T_{\tau_2-\tau_1}f^{(q)})\right)\right\|_2\,ds
  \lesssim \|f^{(q)}-e^{\lambda_1(\tau_2-\tau_1)/2}T_{\tau_2-\tau_1}f^{(q)}\|_4.
\end{equation}
Combining \eqref{4.21} and \eqref{4.22} we get that
\begin{equation}\label{4.24}
  \Sigma(\tau_1,\tau_2)
  \lesssim \|f^{(q)}-e^{\lambda_1(\tau_2-\tau_1)/2}T_{\tau_2-\tau_1}f^{(q)}\|_2
  +\|f^{(q)}-e^{\lambda_1(\tau_2-\tau_1)/2}T_{\tau_2-\tau_1}f^{(q)}\|_4.
\end{equation}
It follows from Fubini's theorem that, for $p=2,4$,
\begin{eqnarray*}
  &&\|U_qf-e^{\lambda_1(\tau_2-\tau_1)/2}T_{\tau_2-\tau_1}U_qf\|_p=
\left\|\int_{0}^\infty e^{-qu}T_uf\,du-e^{(\lambda_1/2+q)(\tau_2-\tau_1)}\int_{\tau_2-\tau_1}^\infty e^{-qu}T_uf\,du\right\|_p\\
   &\le &  \left\|\int_{0}^{\tau_2-\tau_1} e^{-qu}T_uf\,du\right\|_p
   +\left(e^{(\lambda_1/2+q)(\tau_2-\tau_1)}-1\right)\left\|\int_{\tau_2-\tau_1}^\infty e^{-qu}T_uf\,du\right\|_p\\
   &\le &\int_{0}^{\tau_2-\tau_1} e^{-qu}\|T_uf\|_p\,du
   +(e^{(\lambda_1/2+q)(\tau_2-\tau_1)}-1)\int_{\tau_2-\tau_1}^\infty e^{-qu}\|T_uf\|_p\,du.
\end{eqnarray*}
Since $\|T_uf\|_p\le e^{Ku}\|f\|_p$ and $q>K$, we have
\begin{equation}\label{4.16}
  \int_{0}^{\tau_2-\tau_1} e^{-qu}\|T_uf\|_p\,du\le \int_{0}^{\tau_2-\tau_1} e^{-qu}e^{Ku}\,du\|f\|_p\le (\tau_2-\tau_1)\|f\|_p.
\end{equation}
If $\tau_2-\tau_1>t_0$, by \eqref{1.36}, for $u>\tau_2-\tau_1$,
we have $\|T_uf\|_p\lesssim e^{-\lambda_{\gamma(f)}u}\|f\|_2\|a_{t_0}^{1/2}\|_p.$
Thus,
\begin{eqnarray}\label{4.17}
  &&\left(e^{(\lambda_1/2+q)(\tau_2-\tau_1)}-1\right)\int_{\tau_2-\tau_1}^\infty e^{-qu}\|T_uf\|_p\,du\lesssim e^{(\lambda_1/2+q)(\tau_2-\tau_1)}\int_{\tau_2-\tau_1}^\infty e^{-qu}e^{-\lambda_{\gamma(f)}u}\,du\|f\|_2\nonumber\\
  &\lesssim& e^{(\lambda_1/2-\lambda_{\gamma(f)})(\tau_2-\tau_1)}\lesssim \tau_2-\tau_1.
\end{eqnarray}
If
$\tau_2-\tau_1\le t_0$,
then $e^{(\lambda_1/2+q)(\tau_2-\tau_1)}-1\lesssim \tau_2-\tau_1.$
Thus,
\begin{equation}\label{4.18}
  \left(e^{(\lambda_1/2+q)(\tau_2-\tau_1)}-1\right)\int_{\tau_2-\tau_1}^\infty e^{-qu}\|T_uf\|_p\,du\le (e^{(\lambda_1/2+q)(\tau_2-\tau_1)}-1)\|f\|_p\int_{0}^\infty e^{-qu}e^{Ku}\,du\lesssim \tau_2-\tau_1.
\end{equation}
Now, combining \eqref{4.16}--\eqref{4.18}, we obtain that, for $p=2,4$,
$$
\left\|U_qf-e^{\lambda_1(\tau_2-\tau_1)/2}T_{\tau_2-\tau_1}U_qf\right\|_p\lesssim \tau_2-\tau_1.
$$
Now, by \eqref{4.24}, we have
\begin{equation}\label{4.23}
  \Sigma(\tau_1,\tau_2)\le C(\tau_2-\tau_1).
\end{equation}
Thus, by \eqref{4.15} and \eqref{4.23}, we get
\begin{equation}\label{4.4}
  E|G^{1,U_qf}(\tau_2)-G^{1,U_qf}(\tau_1)|^4\le C(\tau_2-\tau_1).
\end{equation}

(2)
We claim that
\begin{equation}\label{4.5}
  E|G^{3,g}(\tau_2)-G^{3,g}(\tau_1)|^4\le C(\tau_2-\tau_1),
\end{equation}
where $C$ is a constant.
To prove \eqref{4.5}, using the same argument as that of leading to \eqref{4.15}, it suffices to show that, for $0\le \tau_1\le \tau_2$,
\begin{equation}\label{5.41}
  E(G^{3,g}(\tau_2)-G^{3,g}(\tau_1))^2\le C(\tau_2-\tau_1).
\end{equation}
Note that
\begin{eqnarray*}
  &&E(G^{3,g}(\tau_2)-G^{3,g}(\tau_1))^2=2\beta_{g,0}-2\beta_{g,\tau_2-\tau_1}\\
  &=&
  2\int_0^\infty e^{-\lambda_1s}(A(I_sg)^2,\phi_1)_m\,ds
   -2e^{-\lambda_1(\tau_2-\tau_1)/2}\int_0^\infty e^{-\lambda_1s}(A(I_sg)(I_{s+\tau_2-\tau_1}g),\phi_1)_m\,ds\\
  &=&2\int_0^\infty e^{-\lambda_1s}
  \left( A(I_sg)(I_sg-e^{-\lambda_1(\tau_2-\tau_1)/2}I_{s+\tau_2-\tau_1}g),\phi_1\right)_m\,ds.
\end{eqnarray*}
By \eqref{1.37}, we have that for any $x\in E$,
$$
 |I_sg(x)|\le \sum_{k: 2\lambda_k<\lambda_1}\sum_{j=1}^{n_k}e^{\lambda_ks}|b_j^k||\phi_j^k(x)|
\lesssim e^{\lambda_{k_0} s}a_{2t_0}(x)^{1/2},
$$
where $k_0=\sup\{k:2\lambda_k<\lambda_1\}$.
By the definition of $I_ug$,
\begin{eqnarray*}
  &&\left|I_sg-e^{-\lambda_1(\tau_2-\tau_1)/2}I_{s+\tau_2-\tau_1}g\right|
  = \left|\sum_{k: 2\lambda_k<\lambda_1}\sum_{j=1}^{n_k}e^{\lambda_ks}(1-e^{(\lambda_k-\lambda_1/2)(\tau_2-\tau_1)})b_j^k\phi_j^{(k)}(x)\right|\\
  &\le &(-\lambda_1/2)(\tau_2-\tau_1)\sum_{k: 2\lambda_k<\lambda_1}\sum_{j=1}^{n_k}e^{\lambda_ks}|b_j^k||\phi_j^{(k)}(x)|
  \lesssim (-\lambda_1/2)(\tau_2-\tau_1)e^{\lambda_{k_0} s}a_{2t_0}(x)^{1/2}.
\end{eqnarray*}
It follows that
\begin{eqnarray*}
  &&E(G^{3,g}(\tau_2)-G^{3,g}(\tau_1))^2\lesssim (-\lambda_1)K(\tau_2-\tau_1)\int_0^\infty e^{-\lambda_1s}e^{2\lambda_{k_0} s}
 ( a_{2t_0},\phi_1)_m\,ds\\
  &=&(-\lambda_1)K(\lambda_1-2\lambda_{k_0})^{-1}(a_{2t_0},\phi_1)_m(\tau_2-\tau_1).
\end{eqnarray*}
Now the proof is complete.
\hfill$\Box$

By Lemma \ref{fin-dim} and Lemma \ref{contin}, we get the following Corollary immediately.
\begin{cor}\label{lem4}
Let $f\in \C_s$, $h\in\C_c$, $g\in\C_l$ and $\mu\in\mathcal{M}_F(E)$.
Suppose that $Y^{1,f}_t$, $Y^{2,h}_t$, $Y^{3,g}_t$, $G^{1,U^qf}$, $G^{2,h}$ and $G^{3,g}$ are defined as in Theorem \ref{Thm1}.
Then, under $\P_{\mu}$, as $t\to\infty$,
\begin{equation}\label{2.20}
   \big(W_t, Y^{1,U_qf}_t,Y^{2,h}_t,Y^{3,g}_t\big)\stackrel{\mathcal{L}(\R_+)}{\longrightarrow} \big(W_\infty,\sqrt{W_\infty}G^{1,U_qf},\sqrt{W_\infty}G^{2,h}\sqrt{W_\infty}G^{3,g}\big).
\end{equation}
\end{cor}

\subsection{The tightness of $\big(W_t, Y^{1,U_qf}_t,Y^{2,h}_t,Y^{3,g}_t\big)_{t>0}$ in $\D(\R^4)$}

Recall that a sequence $(X^n)$ of cadlag processes is called $C$-tight if it is tight,
and if all its weakly convergent limit points are continuous processes.
In this subsection, we will show that $\big(W_t, Y^{1,U_qf}_t,Y^{2,h}_t,Y^{3,g}_t\big)_{t>0}$ is $C$-tight in  $\D(\R^4)$ (with $W_t$, for each $t>0$, being considered as a constant process).
By \cite[Chapter VI, Corollary 3.33]{J.J.}, it suffices to show that $\left(Y^{1,U_qf}_t\right)_{t>0}$, $\left(Y^{2,h}_t\right)_{t>0}$ and $\left(Y^{3,g}_t\right)_{t>0}$ are $C$-tight in  $\D(\R)$.

\subsubsection{The tightness of $\left(Y^{1,U_qf}_t\right)_{t>0}$ in $\D(\R)$}

The main purpose of this subsection is to prove
that $\left(Y^{1,U_qf}_t(\cdot)\right)_{t>0}$ is $C$-tight in  $\D(\R)$.
The next lemma  gives a sufficient condition for the tightness of a sequence $(X^n)_{n\ge 1}$ in $\D(\R^d)$.
\begin{lemma}\label{lem6}
Assume  $(X^n)_{n\ge 1}$ is a sequence of $\D(\R^d)$-valued random variables,
each $X^n$ being defined on the space $\left(\Omega^n,\mathcal{F}^n, \{\mathcal{F}^n_t\}_{t\ge 0}, P^n\right)$.
If $(X^n)$ satisfies the following two conditions:
\begin{enumerate}
  \item [(1)] For all $N>0$,
  \begin{equation}\label{cond1}
    \limsup_{n\to\infty}P^n\left(\sup_{t\le N}|X^n_t|\right)<\infty.
  \end{equation}
  \item [(2)]
  For all $N>0$,
  \begin{equation}\label{cond2}
     \lim_{\theta\to 0}\limsup_{n}\sup_{S,T\in\mathcal{T}_N^n:S\le T\le S+\theta}P^n\left(|X_T^n-X_S^n|\right)=0,
  \end{equation}
  where $\mathcal{T}_N^n$ denotes the set of all $\{\mathcal{F}^n_t\}$-stopping times that are bounded by $N$.
\end{enumerate}
Then, the sequence $(X^n)$ is tight in $\D(\R^d)$.
\end{lemma}
\textbf{Proof:} This follows immediately from Theorem 4.5 in \cite[Chapter VI]{J.J.}.
\hfill$\Box$

 \bigskip

To prove the tightness of $(Y^{1,U_qf}_t(\cdot))_{t>0}$ in $\D(\R)$,
we will check that $Y^{1,U_qf}_t$ satisfies the two conditions above.

\begin{lemma}\label{lem:U1}
If $f\in \C_s$ and $\mu\in\mathcal{M}_F(E)$, then
for any $N>0$,
\begin{equation}\label{5.19}
  \sup_{t>3t_0}\P_{\mu}\left(\sup_{\tau\le N}|Y^{1,U_qf}_t(\tau)|\right)<\infty.
\end{equation}
\end{lemma}
\textbf{Proof:}
In this proof, we always assume that $t>3t_0$. By \eqref{version}, for any $t>0$,
$$\P_{\mu}\left(Y^{1,U_qf}_t(\tau)
=e^{\lambda_1(t+\tau)/2}J_1^f(t+\tau)+e^{(q+\lambda_1/2)(t+\tau)}J_{2}^f(t+\tau),
\forall \tau\ge0\right)=1.$$
First, we consider $J_1^f(t+\tau)$.
Recall that $J_1^f(t)=\langle T_tg,\mu\rangle$, $t\ge 0$.
By \eqref{1.36}, we have
\begin{eqnarray}\label{5.18}
 &&\sup_{\tau\le N}e^{\lambda_1(t+\tau)/2}|J_1^f(t+\tau)|\le  \sup_{\tau\le N}e^{\lambda_1(t+\tau)/2}\langle |T_{t+\tau}g|,\mu\rangle\nonumber\\
 &\lesssim& \sup_{\tau\le N}e^{\lambda_1(t+\tau)/2}e^{-\lambda_{\gamma(g)}(t+\tau)}\|g\|_2\langle a_{t_0}^{1/2},\mu\rangle\nonumber\\
  &\lesssim & e^{(\lambda_1/2-\lambda_{\gamma(g)})t}\|g\|_2
 \lesssim   e^{(\lambda_1/2-\lambda_{\gamma(f)})t}\|f\|_2.
 \end{eqnarray}
Next, we deal with $J_{2}^f(t+\tau)$.
Recall that
$$
  J_{2}^f(t+\tau)=\int_{t+\tau}^\infty e^{-qu}M^{(u)}_{t+\tau}\,du.
$$
Using \eqref{4.10} with $t_1=t$, we have, for $t>3t_0$,
\begin{eqnarray}
  &&\P_{\mu}\left(\sup_{\tau\le N}|J_{2}^f(t+\tau)|\right)
 \le \P_{\mu}\int_t^{\infty}e^{-qu}
  \sup_{\tau\le N}\left(\left|M^{(u)}_{t+\tau}\right|\textbf{1}_{ t+\tau<u}\right)\,du\nonumber\\
  &\le& 2\int_t^{\infty}e^{-qu}
  \sqrt{\int_E \V ar_{\delta_x}\langle f,X_u\rangle\,\mu(dx)}\,du\nonumber\\
  &\lesssim&\int_t^{\infty}e^{-qu}e^{-\lambda_1u/2}\,du\sqrt{\int_E a_{t_0}(x)^{1/2}\,\mu(dx)}\nonumber\\
  &=&(q+\lambda_1/2)^{-1}e^{-(q+\lambda_1/2)t}\sqrt{\int_E a_{t_0}(x)^{1/2}\,\mu(dx)},
  \label{5.7}
\end{eqnarray}
where in the third inequality we use \eqref{2.10}.
It follows that,
$$\sup_{t>3t_0}\P_{\mu}\left(\sup_{\tau\le N}e^{(q+\lambda_1/2)(t+\tau)}J_{2}^f(t+\tau)\right)\le \sup_{t>3t_0}e^{(q+\lambda_1/2)(t+N)}\P_{\mu}\left(\sup_{\tau\le N}|J_{2}^f(t+\tau)|\right)<\infty.$$
The proof is now complete.

\hfill$\Box$

Next, we prove that
\begin{lemma}\label{lem:U2}
If $f\in \C_s$ and $\mu\in\mathcal{M}_F(E)$, then
\begin{equation}\label{6.1}
  \lim_{\theta\to0}\limsup_{t\to\infty} \sup_{S,T\in   \mathcal{T}^t_{N}:S<T<S+\theta}
  \P_{\mu}\left(|Y^{1,U_qf}_t(T)-Y^{1,U_qf}_t(S)|\right)=0,
\end{equation}
where $\mathcal{T}^t_{N}$ is the set of all
$\{\mathcal{G}_{t+\tau}:\tau\ge 0\}$-stoping times that are bounded by $N$.
\end{lemma}
\textbf{Proof:}
In this proof, we always assume that $t>3t_0$.
By \eqref{version}, we have, $\P_{\mu}$-a.s.,
\begin{eqnarray*}
  |Y^{1,U_qf}_t(T)-Y^{1,U_qf}_t(S)|&\le& |e^{\lambda_1(t+T)/2}J_1^f(t+T)-e^{\lambda_1(t+S)/2}J_1^f(t+S)|\\
  &&+|e^{(q+\lambda_1/2)(t+T)}J_{2}^f(t+T)-e^{(q+\lambda_1/2)(t+S)}J_{2}^f(t+S)|\\
   &:=&J_{3,1}(t,T,S)+  J_{3,2}(t,T,S).
\end{eqnarray*}
For $J_{3,1}(t,T,S)$, by \eqref{5.18}, we have that, as $t\to\infty$,
\begin{eqnarray}\label{6.2}
  \P_{\mu}J_{3,1}(t,T,S) &\le & 2 \P_{\mu}\left(\sup_{\tau\le N}e^{\lambda_1(t+\tau)/2}|J_1^f(t+\tau)|\right)
  \lesssim e^{(\lambda_1/2-\lambda_{\gamma(f)})t}\|f\|_2\to 0.
\end{eqnarray}
Note that
\begin{eqnarray*}
  &&J_{3,2}(t,T,S)\\
  &\le& e^{(q+\lambda_1/2)(t+S)}|J_{2}^f(t+T)-J_{2}^f(t+S)|+|e^{(q+\lambda_1/2)(t+T)}-e^{(q+\lambda_1/2)(t+S)}||J_{2}^f(t+T)|\\
  &\le&  e^{(q+\lambda_1/2)(t+N)}|J_{2}^f(t+T)-J_{2}^f(t+S)|+e^{(q+\lambda_1/2)(t+N)}|e^{(q+\lambda_1/2)\theta}-1||J_{2}^f(t+T)|.
\end{eqnarray*}
By \eqref{5.7}, we get that, for $t>3t_0$,
\begin{eqnarray}\label{6.3}
  &&\sup_{S,T\in \mathcal{T}_{N}^t: S<T<S+\theta}e^{(q+\lambda_1/2)(t+N)}
  |e^{(q+\lambda_1/2)\theta}-1|\P_{\mu}|J_{2}^f(t+T)|\nonumber\\
  &\lesssim& e^{(q+\lambda_1/2)(t+N)}|e^{(q+\lambda_1/2)\theta}-1|\P_{\mu}\left(\sup_{\tau\le N}|J_{2}^f(t+\tau)|\right)\nonumber\\
  &\lesssim &|e^{(q+\lambda_1/2)\theta}-1|\to 0, \mbox{ as } \theta\to 0.
\end{eqnarray}
By \eqref{6.2} and \eqref{6.3}, to prove \eqref{6.1}, it suffices to show that
\begin{equation}\label{6.4}
  \lim_{\theta\to0}\limsup_{t\to\infty} \sup_{S,T\in   \mathcal{T}^t_{N}:S<T<S+\theta}
  e^{(q+\lambda_1/2)t}\P_{\mu}|J_{2}^f(t+T)-J_{2}^f(t+S)|=0.
\end{equation}

By the definition of $J_{2}^f$, we have
\begin{eqnarray*}
  &&|J_{2}^f(t+T)-J_{2}^f(t+S)|=\left|\int_{t+T}^\infty e^{-qu}M_{t+T}^{(u)}\,du-\int_{t+S}^\infty e^{-qu}M_{t+S}^{(u)}\,du\right|\\
  &\le&\int_{t+T}^\infty e^{-qu}\left|M_{t+T}^{(u)}-M_{t+S}^{(u)}\right|\,du+\int_{t+S}^{t+T} e^{-qu}|M_{t+S}^{(u)}|\,du\\
  &\le&\int_{t}^\infty e^{-qu}\left|M_{(t+T)\wedge u}^{(u)}-M_{(t+S)\wedge u}^{(u)}\right|\,du+\int_{t+S}^{t+T} e^{-qu}|M_{t+S}^{(u)}|\,du\\
  &:=&J_{4}(t,T,S)+J_{5}(t,T,S).
\end{eqnarray*}
First, we deal with $J_{4}$.
Since $T,S\in \mathcal{T}^t_{N}$, $(t+T)\wedge u$ and $(t+S)\wedge u$ are both
$\{\mathcal{G}_\tau:\tau\ge 0\}$-stopping times.
Thus, by \eqref{5.22}, we have
\begin{eqnarray*}
   &&\P_{\mu} J_{4}(t,T,S)\le \int_{t}^\infty e^{-qu}\sqrt{\P_{\mu}\left|M_{(t+T)\wedge u}^{(u)}-M_{(t+S)\wedge u}^{(u)}\right|^2}\,du\\
   &=&\int_{t}^\infty e^{-qu}\sqrt{\P_{\mu} \left(\langle M^{(u)}\rangle_{(t+T)\wedge u}-\langle M^{(u)}\rangle_{(t+S)\wedge u}\right)}\,du\\
   &=&\int_{t}^\infty e^{-qu}\sqrt{\P_{\mu}\int_{(t+S)\wedge u}^{(t+T)\wedge u}\langle A(T_{u-s}f)^2,X_s\rangle\,ds}\,du\\
   &=&\int_{0}^\infty e^{-q(u+t)}\sqrt{\P_{\mu}\int_{S\wedge u}^{T\wedge u}\langle A(T_{u-s}f)^2,X_{s+t}\rangle\,ds}\,du\\
   &\le&\int_{0}^\infty e^{-q(u+t)}\sqrt{\int_{0}^{N}e^{-\lambda_1(t+s)}\P_{\mu}\left|e^{\lambda_1(t+s)}
   \langle A(T_{u-s}f)^2,X_{s+t}\rangle-( A(T_{u-s}f)^2,\phi_1)_m W_\infty\right|\,ds}\,du\\
   &&\quad+\int_{0}^\infty e^{-q(u+t)}\sqrt{\P_{\mu}\int_{S\wedge u}^{T\wedge u}e^{-\lambda_1(t+s)}
   ( A(T_{u-s}f)^2,\phi_1)_m W_\infty\,ds}\,du\\
   &:=&J_{4,1}(t)+J_{4,2}(t,T,S).
\end{eqnarray*}
Now we consider $J_{4,1}$. Let
$V(u-s,t+s):=\P_{\mu}\left|e^{\lambda_1(t+s)}\langle A(T_{u-s}f)^2,X_{s+t}\rangle-
(A(T_{u-s}f)^2,\phi_1)_m W_\infty\right|$.
Then,
\begin{equation}\label{6.5}
  J_{4,1}(t)\le e^{-(q+\lambda_1/2)t}e^{-\lambda_1N/2}\int_0^\infty e^{-qu}\sqrt{\int_0^N V(u-s,t+s)\,ds}\,du.
\end{equation}
Since $(T_{u-s}f)^2(x)\le e^{K(u-s)}T_{u-s}(f^2)(x)$, we get that, for $t>3t_0$,
\begin{eqnarray*}
  V(u-s,t+s)&\le& e^{\lambda_1(t+s)}\int_ET_{t+s}[A(T_{u-s}f)^2](x)\mu(dx)+K\|(T_{u-s}f)^2\|_2\P_{\mu}( W_\infty )\\
  &\le&e^{\lambda_1(t+s)}e^{K(u-s)}K\int_ET_{t+u}(f^2)(x)\mu(dx)+K\|T_{u-s}f\|_4^2\P_{\mu} (W_\infty)\\
  &\lesssim &e^{\lambda_1(t+s)}e^{K(u-s)}e^{-\lambda_1(t+u)}K\int_Ea_{t_0}(x)^{1/2}\mu(dx)+Ke^{2K(u-s)}\|f\|_4^2\P_{\mu} (W_\infty)\\
  &\lesssim & e^{(K-\lambda_1)(u-s)}+e^{2K(u-s)}\le e^{(K-\lambda_1)u}+e^{2Ks},
\end{eqnarray*}
where in the third inequality we used \eqref{1.36} and the fact that $\|T_{u-s}\|_4\le e^{K(u-s)}$.
Note that
$$\int_0^\infty e^{-qu}\sqrt{\int_0^Ne^{(K-\lambda_1)u}+e^{2Ku}\,ds}\,du\le N^{1/2}\int_0^\infty e^{-(q-K/2+\lambda_1/2)u}+e^{-(q-K)u}\,du<\infty.$$
By Lemma \ref{lem3},  we get that $V(u-s,t+s)\to 0$ as $t\to\infty$. By the dominated convergence theorem, we get that
$$\lim_{t\to\infty}\int_0^\infty e^{-qu}\sqrt{\int_0^N V(u-s,t+s)\,ds}\,du=0.$$
It follows from \eqref{6.5} that
\begin{equation}\label{6.6}
     \lim_{t\to\infty} e^{(q+\lambda_1/2)t}J_{4,1}(t)=0.
\end{equation}
For $J_{4,2}(t,T,S)$, since
$( A(T_{u-s}f)^2,\phi_1)_m\le \|A(T_{u-s}f)^2\|_2\le Ke^{2K(u-s)}\|f\|_4^2\le Ke^{2Ku}\|f\|_4^2$, we have
\begin{eqnarray*}
  J_{4,2}(t,T,S)&\le& \|f\|_4e^{-(q+\lambda_1/2)t}e^{-\lambda_1N/2}\int_{0}^\infty e^{-(q-K)u}\sqrt{\P_{\mu}\left(K(T\wedge u-S\wedge u)W_\infty\right)}\,du \\
  &\lesssim&\theta^{1/2}e^{-(q+\lambda_1/2)t}\int_{0}^\infty e^{-(q-K)u}\,du=(q-K)^{-1}\theta^{1/2}e^{-(q+\lambda_1/2)t},
\end{eqnarray*}
where in the second inequality we used the fact that $T\wedge u-S\wedge u<\theta.$
Thus, we get
\begin{equation}\label{6.7}
  \lim_{\theta\to0}\limsup_{t\to\infty}\sup_{S,T\in   \mathcal{T}^t_{N}:S<T<S+\theta}e^{(q+\lambda_1/2)t}J_{4,2}(t,T,S)=0.
\end{equation}
Combining \eqref{6.6} and \eqref{6.7}, we get
\begin{equation}\label{6.8}
  \lim_{\theta\to0}\limsup_{t\to\infty} \sup_{S,T\in   \mathcal{T}^t_{N}:S<T<S+\theta}
  e^{(q+\lambda_1/2)t}\P_{\mu}J_4(t,T,S)=0.
\end{equation}
Finally, we consider $J_5(t,T,S)$.
By H\"{o}lder's inequality, we get
\begin{eqnarray*}
  \P_{\mu}J_5(t,T,S)&=&\P_{\mu}\int_{t+S}^{t+T} e^{-qu}|M_{t+S}^{(u)}|\,du
  \le  \sqrt{\P_{\mu}\int_{t+S}^{t+T} e^{-2qu}|M_{t+S}^{(u)}|^2\,du}\sqrt{\P_{\mu}(T-S)}\\
  &\le& \theta^{1/2}\sqrt{\int_{t}^{t+N} e^{-2qu}\P_{\mu}|M_{(t+S)\wedge u}^{(u)}|^2\,du}=\theta^{1/2}\sqrt{\int_{t}^{t+N} e^{-2qu}\P_{\mu}\langle M^{(u)}\rangle_{(t+S)\wedge u}\,du}\\
  &\le&\theta^{1/2}\sqrt{\int_{t}^{t+N} e^{-2qu}\P_{\mu}\langle M^{(u)}\rangle_{u}\,du}=\theta^{1/2}\sqrt{\int_{t}^{t+N} e^{-2qu}\int_E\V ar_{\delta_x}\langle f,X_u\rangle\,\mu(dx)\,du}\\
  &\lesssim&\theta^{1/2}\sqrt{\int_{t}^{t+N} e^{-2qu}e^{-\lambda_1u}\,du\int_Ea_{t_0}(x)^{1/2}\,\mu(dx)}\lesssim \theta^{1/2}e^{-(q+\lambda_1/2)t},
\end{eqnarray*}
where in the second to the last inequality we used \eqref{2.10}.
Thus, we get that
\begin{equation}\label{6.9}
  \lim_{\theta\to0}\limsup_{t\to\infty} \sup_{S,T\in   \mathcal{T}^t_{N}:S<T<S+\theta}
  e^{(q+\lambda_1/2)t}\P_{\mu}J_5(t,T,S)=0.
\end{equation}
Combining \eqref{6.8} and \eqref{6.9}, we get \eqref{6.4} immediately.
The proof is now complete.\hfill$\Box$

\bigskip

\begin{lemma}\label{lem:U}
If $f\in \C_s$ and $\mu\in\mathcal{M}_F(E)$, then, under $\P_\mu$,
the family of processes $\left(Y^{1,U_qf}_t(\cdot)\right)_{t>0}$ is $C$-tight in $\D(\R)$.
\end{lemma}
\textbf{Proof:}
It follows from Lemmas \ref{lem:U1} and  \ref{lem:U2} that
$\left(Y^{1,U_qf}_t(\cdot)\right)_{t>0}$ is tight in $\mathbb{D}(\R)$ under $\P_\mu$.
By Corollary \ref{lem4} and the fact that
$\sqrt{W_\infty}G^{1,U_qf}$ is a continuous process, we obtain that
$\left(Y^{1,U_qf}_t(\cdot)\right)_{t>0}$ is $C$-tight in $\mathbb{D}(\R)$ under $\P_\mu$.
\hfill$\Box$

\bigskip

\subsubsection{The tightness of $\left(Y^{2,h}_t\right)_{t>0}$ in $\D(\R)$}

The next lemma will be used to  prove the tightness of $(Y^{2,h}_t(\cdot))_{t>0}$.
\begin{lemma}\label{lem10}
Suppose that $\{C(\tau),\tau\ge 0\}$ and, for each $t>0$, $\{C_t(\tau),\tau\ge 0\}$
are non-decreasing cadlag processes
defined on the space $(\Omega, \mathcal{F},P)$
such that $C_t(0)=C(0)=0$ and
 for all $\tau\ge 0$,
\begin{equation}\label{8.7}
 \lim_{t\to\infty} C_t(\tau)\to C(\tau) \quad \mbox{ in probability.}
\end{equation}
If $C$ is a continuous process,
then
\begin{equation}\label{8.8}
  \lim_{t\to\infty}\delta(C_t,C)=0 \quad \mbox{ in probability},
\end{equation}
where $\delta$ is the metric compatible with the Skorohod topology defined in \cite[Chapter VI, 1.26]{J.J.}.
Moreover, as
$t\to\infty$,
$$C_t-C\stackrel{d}{\longrightarrow}0,$$
which implies that $(C_t)_{t\ge 0}$ is $C$-tight in $\D(\R)$.
\end{lemma}
\textbf{Proof:}
Let $D$ be the subset of all the positive rational numbers.
For any subsequence $(n_k)$, by a diagonal argument,
we can find a further subsequence $(n'_k)$ and a set $\Omega_0\subset \Omega$ with $P(\Omega_0)=1$
such that for $\tau\in D$ and $\omega\in \Omega_0$,
\begin{equation}\label{6.11}
  \lim_{k\to\infty}C_{n_k'}(\tau)(\omega)= C(\tau)(\omega).
\end{equation}
Thus, by \cite[Chapter VI, Theorem 2.15(c)]{J.J.}, we have, for $\omega\in \Omega_0$,
$$\lim_{k\to\infty}\delta(C_{n'_k}(\omega),C(\omega))=0,$$
which implies \eqref{8.8}.
The remaining assertion follows immediately from \eqref{8.8}.\hfill$\Box$

\begin{lemma}\label{lem9}
If $h\in \mathcal{C}_c$ and $\mu\in\mathcal{M}_F(E)$, then the family of processes
$(Y^{2,h}_t(\cdot))_{t>0}$
is $C$-tight in $\mathbb{D}(\R)$ under $\P_\mu$.
\end{lemma}
\textbf{Proof:}
For $h\in\mathbb{C}_c$, we have $T_th=e^{-\lambda_1t/2}h$.
Thus, by \eqref{5.1}, we get that, for $t\ge0$, $\P_{\mu}-$a.s.
\begin{equation*}
  \langle h,X_t\rangle=e^{-\lambda_1t/2}\langle h,X_0\rangle+e^{-\lambda_1t/2}\int_0^t\int_E e^{\lambda_1s/2}h(x)M(ds,dx).
\end{equation*}
Since both sides of the above equation are cadlag, we have
$$\P_{\mu}\left(\langle h,X_t\rangle=e^{-\lambda_1t/2}\langle h,X_0\rangle+e^{-\lambda_1t/2}\int_0^t\int_E e^{\lambda_1s/2}h(x)M(ds,dx),\forall t>0\right)=1.$$
Thus, we have
\begin{eqnarray*}
  Y^{2,h}_t(\tau)&=&t^{-1/2}\langle h,X_0\rangle+t^{-1/2}\int_0^{t+\tau}\int_E e^{\lambda_1s/2}h(x)M(ds,dx)\\
   &=&Y^{2,h}_t(0)+t^{-1/2}\int_t^{t+\tau}\int_E e^{\lambda_1s/2}h(x)M(ds,dx).
  \end{eqnarray*}
Therefore, $\{Y^{2,h}_t(\tau),\tau\ge0\}$ is a square-integrable martingale with
\begin{equation}\label{7.21}
     \langle Y^{2,h}_t\rangle(\tau)=t^{-1}\int_t^{t+\tau}e^{\lambda_1s}\langle Ah^2,X_s\rangle ds.
\end{equation}
By \eqref{1.36}, we have for $t>t_0$,
$$
t^{-1}\P_{\mu}\left(\int_t^{t+\tau}e^{\lambda_1s}\langle Ah^2,X_s\rangle ds\right)=t^{-1}\int_E\int_t^{t+\tau}e^{\lambda_1s}T_s(Ah^2)(x) ds\,\mu(dx)
\lesssim t^{-1}\tau.
$$
Thus, for any $\tau\ge 0$, as $t\to\infty$,
\begin{equation}\label{7.22}
 \langle Y^{2,h}_t\rangle(\tau)\to 0 \quad \mbox{ in $\P_\mu$-probability.}
\end{equation}
Hence, by Lemma \ref{lem10}, $(\langle Y^{2,h}_t\rangle)_{t>0}$ is $C$-tight in $\mathbb{D}(\R)$ under $\P_\mu$.
Since $Y^{2,g}_t(0)=t^{-1/2}e^{-\lambda_1 t/2}\langle g, X_t\rangle\to {\cal N}(0,\rho^2_g)$
in distribution as $t\to\infty$,
we know that $\{Y^{2,h}_t(0),t\ge 0\}$ is tight in $\R$ under $\P_\mu$.
Therefore, by  \cite[Chapter VI, Theorem 4.13]{J.J.}, we get that $\left(Y^{2,h}_t(\cdot)\right)_{t>0}$
is tight in $\mathbb{D}(\R)$ under $\P_\mu$.
By Corollary \ref{lem4} and the fact that
$\sqrt{W_\infty}G^{2,h}$ is a continuous process, we obtain that
$\left(Y^{2,h}_t(\cdot)\right)_{t>0}$ is $C$-tight in $\mathbb{D}(\R)$ under $\P_\mu$.
The proof is now complete.
\hfill$\Box$

\bigskip

\subsubsection{The tightness of $\left(Y^{3,g}_t\right)_{t>0}$ in $\D(\R)$ }

\begin{lemma}\label{tig:lar}
If $g\in \mathcal{C}_l$ and $\mu\in\mathcal{M}_F(E)$, then the family of processes
$(Y^{3,g}_t(\cdot))_{t>0}$
is $C$-tight in $\mathbb{D}(\R)$ under $\P_\mu$.
\end{lemma}
\textbf{Proof:}
Note that
\begin{eqnarray*}
  Y^{3,g}_t(\tau)&=&\sum_{k:\lambda_1>2\lambda_k}\sum_{j=1}^{n_k}e^{(\lambda_1/2-\lambda_k)(t+\tau)}b^k_j\left(H_{t+\tau}^{k,j}-H_{t}^{k,j}\right)\\
&&+\sum_{k:\lambda_1>2\lambda_k}\sum_{j=1}^{n_k}e^{(\lambda_1/2-\lambda_k)(t+\tau)}b^k_j\left(H_{t}^{k,j}-H_{\infty}^{k,j}\right)\\
&:=&Z^{1}_t(\tau)+Z^{2}_t(\tau).
\end{eqnarray*}
 For $Z^{2}_t(\tau)$, it is known (see \cite{RSZ3}) that under $\P_\mu$
$$e^{(\lambda_1/2-\lambda_k)t}\left(H_{t}^{k,j}-H_{\infty}^{k,j}\right)\stackrel{d}{\longrightarrow} G\sqrt{W_\infty},$$
where $G$ is a normal random variable.
It follows that under $\P_\mu$, as $t\to\infty$,
$$e^{(\lambda_1/2-\lambda_k)(t+\cdot)}b^k_j\left(H_{t}^{k,j}-H_{\infty}^{k,j}\right)\stackrel{d}{\longrightarrow}
b^k_jG\sqrt{W_\infty}e^{(\lambda_1/2-\lambda_k)\cdot}.$$
Thus, $e^{(\lambda_1/2-\lambda_k)(t+\cdot)}b^k_j\left(H_{t}^{k,j}-H_{\infty}^{k,j}\right)$
is $C$-tight in $\D(\R)$ under $\P_\mu$.
By \cite[Corollary 3.33]{J.J.}, $(Z^{2}_t)_{t>0}$ is $C$-tight in $\D(\R)$ under $\P_\mu$.
Thus, to prove $(Y^{3,g}_t)_{t>0}$ is tight in $\D(\R)$ under $\P_\mu$,
it suffices to show that $(Z^{1}_t)_{t>0}$ is tight in $\D(\R)$ under $\P_\mu$.

Since $\{H_{t+\tau}^{k,j}-H_{t}^{k,j}:\tau\ge 0\}$ is a martingale under $\P_\mu$,
using $L_p$ maximum inequality, we get
for $\lambda_1>2\lambda_k$,
\begin{eqnarray*}
 \P_{\mu}\left(\sup_{\tau\le N}e^{(\lambda_1/2-\lambda_k)(t+\tau)}\left|H_{t+\tau}^{k,j}-H_{t}^{k,j}\right|\right)&\le & 2e^{(\lambda_1/2-\lambda_k)(t+N)}\sqrt{\P_{\mu}\left(H_{t+N}^{k,j}-H_{t}^{k,j}\right)^2}.
\end{eqnarray*}
By \eqref{5.1}, we have
\begin{equation}\label{5.31}
H_t^{k,j}=\langle \phi_j^{(k)},\mu\rangle+\int_0^t\int_E e^{\lambda_ks}\phi_j^{(k)}(x)\,M(ds,dx).
\end{equation}
Thus,
\begin{equation}\label{5.32}
  \langle H^{k,j}\rangle_t=\int_0^t e^{2\lambda_ks}\langle A(\phi_j^{(k)})^2,X_s\rangle\,ds.
\end{equation}
Therefore, by \eqref{1.36}, we get that, for $t>t_0$,
$$\P_{\mu}\left(H_{t+N}^{k,j}-H_{t}^{k,j}\right)^2
=\int_E\int_t^{t+N} e^{2\lambda_ks}T_s\left( A(\phi_j^{(k)})^2\right)(x)\,ds\,\mu(dx)\lesssim \int_t^{t+N} e^{2\lambda_ks}e^{-\lambda_1s}\,ds
\lesssim e^{(2\lambda_k-\lambda_1)t}.$$
Hence,
\begin{equation}\label{5.36}
  \sup_{t>t_0}\P_{\mu}\left(\sup_{\tau\le N}e^{(\lambda_1/2-\lambda_k)(t+\tau)}\left|H_{t+\tau}^{k,j}-H_{t}^{k,j}\right|\right)<\infty.
\end{equation}
It follows that
\begin{equation}\label{5.33}
  \sup_{t>t_0}\P_{\mu}\left(\sup_{\tau\le N}\left|Z^{1}_t(\tau)\right| \right)\le\sum_{k:\lambda_1>2\lambda_k}\sum_{j=1}^{n_k}|b^k_j|
  \sup_{t>t_0}\P_{\mu}\left(\sup_{\tau<N}e^{(\lambda_1/2-\lambda_k)(t+\tau)}\left|H_{t+\tau}^{k,j}-H_{t}^{k,j}\right|\right)<\infty.
\end{equation}

 Next we prove that
 \begin{equation}\label{5.34}
   \lim_{\theta\to0}\limsup_{t\to\infty}\sup_{T,S\in \mathcal{T}^t_N:0\le T-S\le \theta}\P_{\mu}\left(\left|Z^{1}_t(T)-Z^{1}_t(S)\right| \right)=0,
 \end{equation}
 where $\mathcal{T}^t_{N}$ is the set of all $\{\mathcal{G}_{t+\tau}:\tau\ge 0\}$-stoping times that are bounded by $N$.
 It suffices to show that, for $\lambda_1>2\lambda_k$,
 \begin{equation}\label{5.37}
   \lim_{\theta\to0}\limsup_{t\to\infty}\sup_{T,S\in \mathcal{T}^t_N:0\le T-S\le \theta}\P_{\mu}\left(\left|e^{(\lambda_1/2-\lambda_k)(t+T)}(H^{k,j}_{t+T}-H_t^{k,j})-e^{(\lambda_1/2-\lambda_k)(t+S)}(H^{k,j}_{t+S}-H_t^{k,j})\right| \right)=0.
 \end{equation}
 We note that
 \begin{eqnarray*}
   &&\left|e^{(\lambda_1/2-\lambda_k)(t+T)}(H^{k,j}_{t+T}-H_t^{k,j})-e^{(\lambda_1/2-\lambda_k)(t+S)}(H^{k,j}_{t+S}-H_t^{k,j})\right|  \\
   &\le&  e^{(\lambda_1/2-\lambda_k)(t+S)}\left|H^{k,j}_{t+T}-H_{t+S}^{k,j}\right|
   +e^{(\lambda_1/2-\lambda_k)(t+S)}
     (e^{(\lambda_1/2-\lambda_k)\theta}-1)\left|H^{k,j}_{t+T}-H_t^{k,j}\right|\\
   &\le& e^{(\lambda_1/2-\lambda_k)(t+N)}\left|H^{k,j}_{t+T}-H_{t+S}^{k,j}\right|+
      e^{(\lambda_1/2-\lambda_k)(t+N)}(e^{(\lambda_1/2-\lambda_k)\theta}-1)\sup_{\tau<N}\left|H^{k,j}_{t+\tau}-H_t^{k,j}\right|.
 \end{eqnarray*}
By \eqref{5.36}, we get that, for $t>t_0$,
\begin{equation}\label{5.39}
  e^{(\lambda_1/2-\lambda_k)(t+N)}(e^{(\lambda_1/2-\lambda_k)\theta}-1)\P_{\mu}\left(\sup_{\tau<N}\left|H^{k,j}_{t+\tau}-H_t^{k,j}\right|\right)
  \lesssim e^{(\lambda_1/2-\lambda_k)\theta}-1\to0,
\end{equation}
as $\theta\to0$.
By \eqref{5.32}, we have
\begin{eqnarray*}
&&e^{(\lambda_1/2-\lambda_k)(t+N)}\P_{\mu}\left|H^{k,j}_{t+T}-H_{t+S}^{k,j}\right|
\le e^{(\lambda_1/2-\lambda_k)(t+N)}\sqrt{\P_{\mu}\left|H^{k,j}_{t+T}-H_{t+S}^{k,j}\right|^2} \nonumber  \\
& =&  e^{(\lambda_1/2-\lambda_k)(t+N)}\sqrt{\P_{\mu}\left(\langle H^{k,j}\rangle_{t+T}-\langle H^{k,j}\rangle_{t+S}\right)}\nonumber\\
&=&e^{(\lambda_1/2-\lambda_k)(t+N)}\sqrt{\P_{\mu}\int_{t+S}^{t+T} e^{2\lambda_ks}\langle A(\phi_j^{(k)})^2,X_s\rangle\,ds}\nonumber\\
&\lesssim & \sqrt{\P_{\mu}\int_{t+S}^{t+T} e^{\lambda_1s}\langle A(\phi_j^{(k)})^2,X_s\rangle\,ds}\nonumber\\
&\le & \sqrt{\int_{t}^{t+N}\P_{\mu}\left| e^{\lambda_1s}\langle A(\phi_j^{(k)})^2,X_s\rangle-
(A(\phi_j^{(k)})^2,\phi_1 )_mW_\infty\right|\,ds
+\theta( A(\phi_j^{(k)})^2,\phi_1 )_m\P_{\mu}(W_\infty) }.
\end{eqnarray*}
By Lemma \ref{lem3},
$$\lim_{t\to\infty}\int_{t}^{t+N}\P_{\mu}\left| e^{\lambda_1s}\langle A(\phi_j^{(k)})^2,X_s\rangle-
(A(\phi_j^{(k)})^2,\phi_1 )_m W_\infty\right|\,ds =0.$$
Thus,
\begin{eqnarray}\label{5.40}
  &&\lim_{\theta\to0}\limsup_{t\to\infty}\sup_{T,S\in \mathcal{T}^t_N:0\le T-S\le \theta}
  e^{(\lambda_1/2-\lambda_k)(t+N)}\P_{\mu}\left|H^{k,j}_{t+T}-H_{t+S}^{k,j}\right|\nonumber\\
  &\lesssim & \lim_{\theta\to0}\sqrt{\theta( A(\phi_j^{(k)})^2,\phi_1 )_m\P_{\mu}(W_\infty) } =0.
\end{eqnarray}
Combining \eqref{5.39} and \eqref{5.40}, we get \eqref{5.37}.

By Corollary \ref{lem4} and the fact that
$\sqrt{W_\infty}G^{3,g}$ is a continuous process, we obtain that
$\left(Y^{3,g}_t(\cdot)\right)_{t>0}$ is $C$-tight in $\mathbb{D}(\R)$ under $\P_\mu$.
The proof is now complete.\hfill$\Box$

\begin{singlespace}

\end{singlespace}

\vskip 0.2truein
\vskip 0.2truein

\noindent{\bf Yan-Xia Ren:} LMAM School of Mathematical Sciences \& Center for
Statistical Science, Peking
University,  Beijing, 100871, P.R. China. Email: {\texttt
yxren@math.pku.edu.cn}

\smallskip
\noindent {\bf Renming Song:} Department of Mathematics,
University of Illinois,
Urbana, IL 61801, U.S.A.
Email: {\texttt rsong@math.uiuc.edu}

\smallskip

\noindent{\bf Rui Zhang:} LMAM School of Mathematical Sciences, Peking
University,  Beijing, 100871, P.R. China. Email: {\texttt
ruizhang8197@gmail.com}

\end{doublespace}
\end{document}